\RequirePackage{ifpdf}
\ifpdf 
\documentclass[pdftex]{sigma}
\else
\documentclass{sigma}
\fi

\usepackage[mathscr]{eucal}
\usepackage{enumerate}
\usepackage{xspace}
\usepackage[thinlines]{easymat}

\numberwithin{equation}{section}

\newtheorem{Theorem}{Theorem}[section]

{\theoremstyle{definition}
\newtheorem{Example}[Theorem]{Example}
\newtheorem{Definition}[Theorem]{Definition}
\newtheorem{Remark}[Theorem]{Remark}
}

\def\&{\wedge}

\newcommand{\alp}{\alpha}

\newcommand{\eps}{\varepsilon}

\newcommand{\calC}{{\mathcal C}}

\newcommand{\calH}{{\mathcal H}}

\newcommand{\scB}{{\EuScript B}}
\newcommand{\scD}{{\EuScript D}}

\newcommand{\scF}{{\EuScript F}}

\newcommand{\scU}{{\EuScript U}}
\newcommand{\scX}{{\EuScript X}}
\newcommand{\scY}{{\EuScript Y}}

\newcommand{\R}{\mathbb{R}}
\newcommand{\C}{\mathbb{C}}

\newcommand{\bb}{\mathbb}
\newcommand{\mbf}{\mathbf}
\newcommand{\ds}{\displaystyle}

\begin{document}


\renewcommand{\thefootnote}{$\star$}

\renewcommand{\PaperNumber}{095}

\FirstPageHeading

\ShortArticleName{Geometry of Control-Af\/f\/ine Systems}

\ArticleName{Geometry of Control-Af\/f\/ine Systems\footnote{This paper is a
contribution to the Special Issue ``\'Elie Cartan and Dif\/ferential Geometry''. The
full collection is available at
\href{http://www.emis.de/journals/SIGMA/Cartan.html}{http://www.emis.de/journals/SIGMA/Cartan.html}}}

\Author{Jeanne N. CLELLAND~$^\dag$, Christopher G. MOSELEY~$^\ddag$ and George R. WILKENS~$^\S$}

\AuthorNameForHeading{J.N. Clelland, C.G. Moseley and G.R. Wilkens}

\Address{$^\dag$~Department of Mathematics, 395 UCB, University of
Colorado, Boulder, CO 80309-0395, USA}
\EmailD{\href{mailto:Jeanne.Clelland@colorado.edu}{Jeanne.Clelland@colorado.edu}}

\Address{$^\ddag$~Department of Mathematics and Statistics, Calvin College, Grand Rapids, MI 49546, USA}
\EmailD{\href{mailto:cgm3@calvin.edu}{cgm3@calvin.edu}}

\Address{$^\S$~Department of Mathematics, University of Hawaii at Manoa, 2565 McCarthy Mall,\\
 \hphantom{$^\S$}~Honolulu, HI 96822-2273, USA}
\EmailD{\href{mailto:grw@math.hawaii.edu}{grw@math.hawaii.edu}}

\ArticleDates{Received April 02, 2009, in f\/inal form September 28, 2009;  Published online October 07, 2009}

\Abstract{Motivated by control-af\/f\/ine systems in optimal control theory, we introduce the notion of a {\em point-affine distribution} on a manifold $\scX$~-- i.e., an af\/f\/ine distribution $\scF$ together with a distinguished vector f\/ield contained in $\scF$.  We compute local invariants for point-af\/f\/ine distributions of constant type when $\dim(\scX)=n$, $\text{rank}(\scF)=n-1$, and when $\dim(\scX)=3$, $\text{rank}(\scF)=1$.
Unlike linear distributions, which are characterized by integer-valued invariants~-- namely, the rank and growth vector~-- when $\dim(\scX)\leq 4$, we f\/ind local invariants depending on arbitrary functions even for rank~1 point-af\/f\/ine distributions on manifolds of dimension~2.}


\Classification{58A30; 53C17; 58A15; 53C10}

\renewcommand{\thefootnote}{\arabic{footnote}}
\setcounter{footnote}{0}

\section{Introduction}

In \cite{CM06} and \cite{CMW07}, we introduced the notion of {\em sub-Finsler geometry} as a natural geometric setting for certain problems in optimal control theory.  In local coordinates, a control system may be represented as an underdetermined system of ordinary dif\/ferential equations of the form
\begin{gather}
   \dot{x} = f(x, u), \label{gencontrol}
\end{gather}
where $x \in \bb{R}^n$ represents the \emph{state vector} of the system and $u \in
\bb{R}^s$ represents the {\em control vector}, i.e., variables which may be
specif\/ied freely in order to ``steer" the system in a desired
direction.  More generally, $x$ and $u$ may take values in an
$n$-dimensional manifold $\scX$ and an $s$-dimensional manifold $\scU$,
respectively.  Typically there are constraints on
how the system may be ``steered" from one state to another, so that
$s < n$.  The
systems of greatest interest are {\em controllable}, i.e., given any
two states $x_1$, $x_2$, there exists a solution curve of
\eqref{gencontrol} connecting $x_1$ to $x_2$.

The control systems for which sub-Finsler (or sub-Riemannian) geometry are most relevant are {\em control-linear} systems. These are systems for which the right-hand side of \eqref{gencontrol} is linear in the control variables $u$ and depends smoothly on the state variables $x$; i.e., systems of the form
\begin{gather}
   \dot{x} = A(x) u, \label{qlincontrol}
\end{gather}
where $A(x)$ is an $n \times s$ matrix whose entries are smooth functions of $x$.  For such a system,
\emph{admissible} paths in the state space are those for which the tangent
vector to the path at each point $x \in \scX$ is contained in the
subspace $\scD_x \subset T_x\scX$ determined by the image of the matrix~$A(x)$.  If this matrix has constant (maximal)
rank $s$ on $\scX$, then $\scD$ is a rank $s$ linear
distribution\footnote{We use the adjective ``linear'' to distinguish the usual notion of a distribution $\scD$, where each f\/iber $\scD_x$ is a~linear subspace of $T_x\scX$, from that of an af\/f\/ine distribution $\scF$.  Unless otherwise specif\/ied, the word ``distribution'' will refer to a~linear distribution throughout this paper.}
 on~$\scX$.  Thus the admissible paths in the state space are precisely the {\em horizontal curves} of the
distribution~$\scD$, i.e., curves whose tangent vectors at each point
are contained in $\scD$.

A sub-Finsler (resp., sub-Riemannian) metric on the distribution $(\scX, \scD)$ is def\/ined by speci\-fying a Finsler (resp., Riemannian) metric on each of the subspaces $\scD_x$.  Such a metric may represent a cost function for the control problem, and geodesics for the metric may be thought of as stationary trajectories for the corresponding {\em optimal control} problem.

While the class \eqref{qlincontrol} of control-linear systems contains many interesting examples, a much broader class of interest is that of {\em control-affine} systems, i.e., systems of the form
\begin{gather}
   \dot{x} = a_0(x) + A(x) u, \label{alincontrol}
\end{gather}
where $A(x)$ is a smoothly-varying $n \times s$ matrix and $a_0(x)$ is a smooth vector f\/ield on $\scX$.  (The vector f\/ield $a_0(x)$ is known as the {\em drift} vector f\/ield.)

\begin{Example}
The class of control-af\/f\/ine systems contains as a sub-class the {\em linear} control systems.  These are systems of the form
\begin{gather}\label{lincontrol}
  \dot{x} = Ax + Bu,
\end{gather}
where $A$ is a constant $n \times n$ matrix and $B$ is a constant $n \times s$ matrix.  Note that despite the terminology, the linear control system \eqref{lincontrol} is not control-linear unless $A=0$.
\end{Example}

The geometric structure corresponding to a control-af\/f\/ine system is an {\em affine distribution} on the state space $\scX$.  More precisely (see \cite{Elkin98, Elkin99}):

\begin{Definition}
A {\em rank $s$ affine distribution} $\scF$ on an $n$-dimensional manifold $\scX$ is a smoothly-varying family of $s$-dimensional, af\/f\/ine linear subspaces $\scF_x \subset T_x\scX$.  We will say that $\scF$ is {\em strictly affine} if none of the af\/f\/ine subspaces $\scF_x \subset T_x\scX$ are linear subspaces.  Associated to an af\/f\/ine distribution $\scF$ is the {\em direction distribution}
\[ L_{\scF} = \{ \xi_1 - \xi_2 \mid \xi_1, \xi_2 \in \scF \}. \]
\end{Definition}
Note that $L_{\scF}$ is a rank $s$ linear distribution on $\scX$.

In this paper, we will limit our attention to distributions (both linear and af\/f\/ine) which are {\em bracket-generating} (or {\em almost bracket-generating} in the af\/f\/ine case) and have {\em constant type}.
In order to make these notions precise, we need the notion of the {\em growth vector} of a distribution (cf.~\cite{Montgomery02}).  Let $\scD$ be a linear distribution on a manifold $\scX$, and let $\scD$ also denote the sheaf of smooth vector f\/ields on $\scX$ which are local sections of $\scD$.  The iterated Lie brackets of vector f\/ields in $\scD$ generate a f\/lag of subsheaves
\[ \scD = \scD^1 \subset \scD^2 \subset \cdots \subset T\scX, \]
def\/ined by $\scD^1 = \scD$, and
\[ \scD^{i+1} = \scD^i + [\scD, \scD^i], \qquad i \geq 1. \]
At each point $x\in \scX$, this f\/lag of subsheaves gives a f\/lag of subspaces of $T_x\scX$:
\[ \scD_x \subset \scD_x^2 \subset \cdots \subset T_x\scX. \]

\begin{Definition}\label{def:d-defs}
Let $\scD^\infty = \cup_{i \geq 1} \scD^i \subset T\scX.$
\begin{itemize}\itemsep=0pt
\item The smallest integer $r = r(x)$ such that $\scD^r_x = \scD^{\infty}_x$ is called the {\em step} of the distribution at~$x$.
\item The distribution $\scD$ is {\em bracket-generating} if $\scD^\infty = T\scX$.
\item Set $n_i(x) = \dim \scD_x^i$.  The {\em growth vector} of $\scD$ at $x$ is the integer list $(n_1(x)$, $n_2(x)$, $\ldots$, $n_r(x))$, where $r(x)$ is the step of $\scD$ at $x$.
\item The distribution $\scD$ has {\em constant type} if $n_i(x)$ is constant on $\scX$ for all $i$; i.e., if the growth vector of $\scD$ is constant on $\scX$.
\end{itemize}
\end{Definition}

\begin{Remark}
When working with distributions of constant type, it is customary to consider only bracket-generating distributions.
For if $\scD$ is not bracket-generating, then $\scD^\infty$ is a Frobenius distribution on $\scX$.  The integral manifolds of $\scD^\infty$ def\/ine a local foliation of $\scX$ with the property that any horizontal curve of $\scD$ is contained in a single leaf of the foliation, and the restriction of $\scD$ to each leaf of the foliation is bracket-generating.
\end{Remark}

We can def\/ine similar notions for af\/f\/ine distributions: let $\scF$ be an af\/f\/ine distribution on a~manifold $\scX$, with direction distribution $L_{\scF}$. Let $\scF$ also denote the sheaf of smooth vector f\/ields on $\scX$ which are local sections of $\scF$. The f\/lag of subsheaves
\[ \scF = \scF^1 \subset \scF^2 \subset \cdots \subset T\scX \]
may be def\/ined in the same manner as for distributions: set $\scF^1 = \scF$, and
\[ \scF^{i+1} = \scF^{i} + [\scF, \scF^i], \qquad i \geq 1. \]
At each point $x \in \scX$, this gives a f\/lag of af\/f\/ine subspaces of $T_x\scX$:
\[ \scF_x \subset \scF_x^2 \subset \cdots \subset T_x\scX. \]
Each af\/f\/ine distribution $\scF^i$ in the f\/lag has an associated direction distribution $L_{\scF^i}$; these distributions clearly have the property that
\[ L_{\scF}^i \subset L_{\scF^i}. \]
The step, the growth vector, and the notion of bracket-generating are def\/ined for af\/f\/ine distributions in the same way as for linear distributions.  But we will also want to consider af\/f\/ine distributions with a slightly weaker bracket-generating property.

\begin{Definition}\label{def:ad-almost-bracket-generating}
An af\/f\/ine distribution $\scF$ on an $n$-dimensional manifold $\scX$ is {\em almost bracket-generating} if $\text{rank}(\scF^\infty) = n-1$, and for each $x \in \scX$ and any $\xi(x) \in \scF_x$,
$\text{span}(\xi(x), (L_{\scF^\infty})_x) = T_x\scX$.
\end{Definition}

We also impose an additional condition in the def\/inition of constant type~-- namely, that $\scF$ is strictly af\/f\/ine, and that each element $\scF^i$ is either strictly af\/f\/ine or a linear distribution.  This condition is ref\/lected in the second condition of the following def\/inition.
\begin{Definition}\label{def:ad-constant-type}
The af\/f\/ine distribution $\scF$ has {\em constant type} if
\begin{itemize}\itemsep=0pt
\item $n_i(x)$ is constant on $\scX$ for all $i$; i.e., the growth vector of $\scF$ is constant on $\scX$, and
\item for any section $\xi$ of $\scF$, $\dim (\text{span}(\xi(x), (L_{\scF^i})_x))$ is constant on $\scX$ for all $i$.
\end{itemize}
\end{Definition}

The well-known theorems of Frobenius, Pfaf\/f, and Engel (cf.\ Section~\ref{classical-results}) imply that linear distributions of constant type on manifolds of dimension $n \leq 4$ have only integer-valued invariants: specif\/ically, any two linear distributions of constant type on manifolds of the same dimension, with the same rank and growth vector, are locally equivalent via a dif\/feomorphism of the underlying manifolds. Beginning with $n=5$ and $s=2$, local invariants depending on arbitrary functions appear: the famous paper of \'E. Cartan \cite{Cartan10} describes local invariants of rank 2 distributions on 5-manifolds with growth vector $(2,3,5)$.  More recently, Bryant \cite{Bryant06} has described local invariants of rank 3 distributions on 6-manifolds with growth vector $(3,6)$, and Doubrov and Zeleneko \cite{DZ07, DZ08} have given a fairly comprehensive treatment of maximally nonholonomic distributions of ranks~2 and~3 on manifolds of arbitrary dimension.

The study of af\/f\/ine distributions is more recent and less extensive. Elkin \cite{Elkin98, Elkin99} has studied equivalence of af\/f\/ine distributions on manifolds of dimension $n \leq 4$ under dif\/feomorphisms of the underlying manifolds, resulting in local normal forms for the associated control systems.  Local invariants appear in lower dimensions than for linear distributions: the f\/irst local invariants arise in the case of rank 1 af\/f\/ine distributions on 3-dimensional manifolds.  Local normal forms for generic rank 1 af\/f\/ine distributions on $\mathbb{R}^3$ have also been obtained by Agrachev \cite{Agrachev98} and Wilkens~\cite{Wilkens99}, and for generic rank 1 af\/f\/ine distributions on $\mathbb{R}^n$ and generic rank 2 af\/f\/ine distributions on $\mathbb{R}^4$ and $\mathbb{R}^5$ by Agrachev and Zelenko \cite{AZ07}.

Given a control-af\/f\/ine system \eqref{alincontrol}, one can canonically associate to it an af\/f\/ine distribution~$\scF$ on the state space $\scX$: if we let $a_1(x), \ldots, a_s(x)$ denote the columns of the matrix $A(x)$, then
\[ \scF_x = \left\{a_0(x) + \sum_{i=1}^s u^i a_i(x) \mid u^1, \ldots, u^s \in \R \right\}. \]
Elkin def\/ines two systems
\begin{gather}
\dot{x}   = a_0(x) + \sum_{i=1}^s u^i a_i(x), \label{control-sys-1} \\
\dot{y}   = b_0(y) + \sum_{i=1}^s v^i b_i(y) \label{control-sys-2}
\end{gather}
on state spaces $\scX$, $\scY$, respectively,
to be {\em affine equivalent} if there exists a dif\/feomorphism $\psi:\scX \to \scY$ such that any absolutely continuous curve $x(t)$ in $\scX$ is a solution of \eqref{control-sys-1} if and only if the curve $y(t) = \psi(x(t))$ is a solution of \eqref{control-sys-2}.  More specif\/ically, such a dif\/feomorphism must have the property that
\begin{gather}
\psi_*(a_0(x))  = b_0(\psi(x)) + \sum_{j=1}^s \lambda^j_0(x) b_j(\psi(x)), \label{drift-transformation-1} \\
\psi_*(a_i(x))   = \sum_{j=1}^s \lambda^j_i(x) b_j(\psi(x)),  \qquad 1 \leq i \leq s\nonumber
\end{gather}
for some functions $\lambda^j_0$, $\lambda^j_i$ on $\scX$.  In terms of the associated af\/f\/ine distributions $\scF_{\scX}$, $\scF_{\scY}$, this is equivalent to the statement that $\psi_*(\scF_{\scX}) = \scF_{\scY}$.

This notion of equivalence is quite natural, especially in terms of the geometry of af\/f\/ine distributions.  But with an eye towards studying {\em optimal} control problems of the form \eqref{alincontrol}, we propose a slightly more restrictive def\/inition of equivalence.  According to \eqref{drift-transformation-1}, the drift vector f\/ield $a_0(x)$ in \eqref{alincontrol} may be replaced by any vector f\/ield of the form $a_0(x) + \sum_{j=1}^s \lambda^j_0(x) a_j(x)$, and the resulting control system will be af\/f\/ine equivalent to the original.  But in practice, there is often a preferred choice for the drift vector f\/ield, corresponding to a zero value for some physical control inputs.  This is particularly true in optimal control, where there is typically a~specif\/ic control input whose cost function is minimal.  In this paper, therefore, we will study the following type of geometric object, corresponding to a control-af\/f\/ine system \eqref{alincontrol} with a~f\/ixed drift vector f\/ield $a_0$:

\begin{Definition}\label{point-affine-distrib-def}
A {\em point-affine distribution} on a manifold $\scX$ is an af\/f\/ine distribution $\scF$ on $\scX$, together with a distinguished vector f\/ield $a_0 \in \scF$.
\end{Definition}

We will consider the equivalence problem for point-af\/f\/ine distributions and their associated control-af\/f\/ine systems with respect to the following notion of equivalence:

\begin{Definition}
The point-af\/f\/ine distributions
\[ \scF_{\scX} = a_0 + \text{span}\left(a_1, \ldots, a_s \right), \qquad \scF_{\scY} = b_0 + \text{span}\left(b_1, \ldots, b_s \right) \]
on the manifolds $\scX, \scY$ (corresponding to the control-af\/f\/ine systems \eqref{control-sys-1} and \eqref{control-sys-2}, respectively) will be called {\em point-affine equivalent} if there exists a dif\/feomorphism $\psi:\scX \to \scY$ such that
\begin{gather*}
\psi_*(a_0(x))   = b_0(\psi(x)), \\ 
\psi_*(a_i(x))   = \sum_{j=1}^s \lambda^j_i(x) b_j(\psi(x)),  \qquad 1 \leq i \leq s
\end{gather*}
for some functions $\lambda^j_i$ on $\scX$.
\end{Definition}
In other words, there must exist a dif\/feomorphism that preserves both the af\/f\/ine distributions {\em and} the distinguished vector f\/ields.

The primary goal of this paper is the local classif\/ication of point-af\/f\/ine distributions of constant type with respect to point-af\/f\/ine equivalence.  We will use Cartan's method of equivalence to compute local invariants for point-af\/f\/ine distributions.  The remainder of the paper is organized as follows: in Section~\ref{review}, we review some well-known results on the geometry of linear distributions of constant type, as well as some of Elkin's results on af\/f\/ine distributions, in order to put our results in context.  In Section~\ref{results} we compute local invariants for point-af\/f\/ine distributions of constant type when $\dim (\scX)=n$, $\text{rank}(\scF)=n-1$, and when $\dim (\scX)=3$, $\text{rank}(\scF)=1$.  In Section~\ref{future} we discuss examples and potential directions for further study.

\section{Review of prior results}\label{review}

\subsection{Classical normal form results for linear distributions}\label{classical-results}

The following results are well-known; see, e.g., \cite{BCG3} or \cite{IL03}.  We include them here primarily to provide historical context for our results; some of them will also be used in the proofs in Section~\ref{results}.

\begin{Definition}
A distribution $\scD$ on a manifold $\scX$ is called {\em completely integrable} (or, more succinctly, {\em integrable}), if $\scD^{\infty} = \scD$.
\end{Definition}

\begin{Theorem}[Frobenius]\label{frob-thm}
 Let $\scD$ be a rank $s$ distribution on an $n$-dimensional manifold $\scX$.  If $\scD$ is completely integrable, then in a sufficiently small neighborhood of any point $x \in \scX$, there exist local coordinates $(x^1, \ldots, x^n)$ such that
\begin{gather*}
 \scD = \text{\rm span}\left(\frac{\partial}{\partial x^1}, \ldots, \frac{\partial}{\partial x^s}\right),
\end{gather*}
or, equivalently,
\begin{gather*}
 \scD^{\perp} = \{dx^{s+1}, \ldots, dx^n\}. 
\end{gather*}
\end{Theorem}

\begin{Theorem}[Pfaf\/f]\label{pfaff-thm}
 Let $\scD$ be a rank $n$ distribution on an $(n+1)$-dimensional manifold $\scX$, and let $\theta$ be a nonvanishing $1$-form on $\scX$ such that $\scD = \{\theta\}^{\perp}$.  Let $k$ be the integer defined by the conditions
\[ \theta \& (d\theta)^k  \neq 0, \qquad \theta \& (d\theta)^{k+1}  = 0. \]
$(k$ is called the {\em Pfaff rank} of $\theta.)$  In a sufficiently small neighborhood of any point $x \in \scX$ on which $k$ is constant, there exist local coordinates $(x^0, \ldots, x^n)$ such that
\begin{gather*}
\theta = \begin{cases}
dx^1 + x^2\, dx^3 + \cdots + x^{2k}\, dx^{2k+1} & \text{if} \ \ (d\theta)^{k+1}=0, \\
 x^0\, dx^1 + x^2\, dx^3 + \cdots + x^{2k}\, dx^{2k+1} & \text{if} \ \ (d\theta)^{k+1} \neq 0.
 \end{cases}
\end{gather*}
\end{Theorem}

\begin{Theorem}[Engel]\label{engel-thm}
 Let $\scD$ be a rank $2$ distribution of constant type on a $4$-dimensional manifold $\scX$, with growth vector $(2,3,4)$.  Then in a sufficiently small neighborhood of any point $x \in \scX$, there exist local coordinates $(x^1, x^2, x^3, x^4)$ such that
\begin{gather*}
\scD^{\perp} = \{dx^2 - x^3\, dx^1, \, dx^3 - x^4\, dx^1\}.
\end{gather*}
\end{Theorem}

\subsection{Local classif\/ication of rank 1 af\/f\/ine systems\\ on 2- and 3-dimensional manifolds under af\/f\/ine equivalence} \label{elkin-results}

Linear distributions of rank 1 have no local invariants, as all nonvanishing vector f\/ields on a~manifold are well-known to be locally equivalent.  Af\/f\/ine distributions, however, may have local invariants even in rank 1.  The following results are due to Elkin~\cite{Elkin98, Elkin99}.

\begin{Theorem}[Elkin]\label{elkin-2d-thm}
 Let $\scF$ be a rank $1$ strictly affine distribution of constant type on a~$2$-dimensional manifold $\scX$.
\begin{enumerate}\itemsep=0pt
\item If $\scF$ is almost bracket-generating, then in a sufficiently small neighborhood of any point $x \in \scX$, there exist local coordinates $(x^1, x^2)$ such that
\[ \scF = \frac{\partial}{\partial x^1} + \text{\rm span}\left(\frac{\partial}{\partial x^2}\right). \]
\item If $\scF$ is bracket-generating, then in a sufficiently small neighborhood of any point $x \in \scX$, there exist local coordinates $(x^1, x^2)$ such that
\[ \scF = x^2\, \frac{\partial}{\partial x^1} + \text{\rm span}\left(\frac{\partial}{\partial x^2}\right). \]
\end{enumerate}
\end{Theorem}
Note that, for reasons of dimension, any strictly af\/f\/ine distribution of rank 1 on a 2-manifold is either almost bracket-generating or bracket-generating.

\begin{Theorem}[Elkin]\label{elkin-3d-thm}
 Let $\scF$ be a rank $1$ strictly affine distribution of constant type on a~$3$-dimensional manifold $\scX$.
\begin{enumerate}\itemsep=0pt
\item\label{elkin-3d-case1} If $\scF$ is almost bracket-generating, then in a sufficiently small neighborhood of any point $x \in \scX$, there exist local coordinates $(x^1, x^2, x^3)$ such that
\[ \scF = \left(\frac{\partial}{\partial x^1} + x^3\, \frac{\partial}{\partial x^2}\right) + \text{\rm span}\left(\frac{\partial}{\partial x^3}\right). \]
\item \label{elkin-3d-case2}If $\scF$ is bracket-generating and $L_{\scF^2}$ $($which must have rank~$2)$ is Frobenius, then in a sufficiently small neighborhood of any point $x \in \scX$, there exist local coordinates $(x^1, x^2, x^3)$ such that
\[ \scF = \left(x^2\, \frac{\partial}{\partial x^1} + x^3\, \frac{\partial}{\partial x^2}\right) + \text{\rm span}\left(\frac{\partial}{\partial x^3}\right). \]
\item\label{elkin-3d-case3} If $\scF$ is bracket-generating and $L_{\scF^2}$ is not Frobenius, then in a sufficiently small neighborhood of any point $x \in \scX$, there exist local coordinates $(x^1, x^2, x^3)$ such that
\[ \scF = \frac{\partial}{\partial x^1} + \text{\rm span}\left(x^3\, \frac{\partial}{\partial x^1} + \frac{\partial}{\partial x^2} + H\, \frac{\partial}{\partial x^3} \right), \]
where $H$ is an arbitrary function on $\scX$ satisfying $\frac{\partial H}{\partial x^1} \neq 0$.
\end{enumerate}
\end{Theorem}
Note that if a rank 1 strictly af\/f\/ine distribution $\scF$ on a 3-manifold is not bracket-generating or almost bracket-generating, then $\text{span}(\scF)$ is a rank 2 Frobenius distribution $\scD_{\scF}$ on $\scX$, and $\scF$ may be regarded (at least locally) as living on the 2-dimensional leaves of the foliation which is tangent to $\scD_\scF$.

\section[Local point-affine equivalence results for point-affine distributions]{Local point-af\/f\/ine equivalence results \\ for point-af\/f\/ine distributions} \label{results}

\subsection{Equivalence problem for rank  1  point-af\/f\/ine distributions on 2-manifolds}\label{2-1-results}

For point-af\/f\/ine equivalence, the f\/irst local invariants appear in the simplest nontrivial case: rank 1 point-af\/f\/ine distributions on 2-dimensional manifolds.

\begin{Theorem}\label{2-1-thm}
Let $\scF$ be a rank $1$ strictly affine point-affine distribution of constant type on a~$2$-dimensional manifold $\scX$.
\begin{enumerate}\itemsep=0pt
\item If $\scF$ is almost bracket-generating, then in a sufficiently small neighborhood of any point $x \in \scX$, there exist local coordinates $(x^1, x^2)$ such that
\[ \scF = \frac{\partial}{\partial x^1} + \text{\rm span}\left(\frac{\partial}{\partial x^2}\right). \]
\item If $\scF$ is bracket-generating, then in a sufficiently small neighborhood of any point $x \in \scX$, there exist local coordinates $(x^1, x^2)$ such that
\[ \scF = x^2 \left( \frac{\partial}{\partial x^1} + J \frac{\partial}{\partial x^2} \right) + \text{\rm span}\left(\frac{\partial}{\partial x^2}\right), \]
where $J$ is an arbitrary function on $\scX$.
\end{enumerate}
\end{Theorem}

\begin{proof}
We will employ Cartan's method of equivalence in this proof.
Let $\scF$ be a rank 1 point-af\/f\/ine distribution of constant type on a 2-dimensional manifold $\scX$.  A local framing~-- i.e., a pair of everywhere linearly independent vector f\/ields $(v_1, v_2)$ on $\scX$~-- will be called {\em admissible} if $v_1$ is the distinguished vector f\/ield and
\[ \scF = v_1 + \text{span}(v_2). \]
If $(v_1, v_2)$ is any admissible framing on $\scX$, then any other admissible framing $(\tilde{v}_1, \tilde{v}_2)$ has the form{\samepage
\begin{gather*}
\tilde{v}_1   = v_1 ,\qquad
\tilde{v}_2   = b_2 v_2 ,
\end{gather*}
for some nonvanishing function $b_2$ on $\scX$.}

A local coframing $(\bar{\eta}^1, \bar{\eta}^2)$ on $\scX$ will be called {\em $0$-adapted} if it is the dual coframing of some admissible framing.
Any two 0-adapted coframings on $\scX$ vary by a transformation of the form
\begin{gather}\label{2-1-g0-group}
\begin{bmatrix} \tilde{\bar{\eta}}^1 \\[0.05in] \tilde{\bar{\eta}}^2 \end{bmatrix} =
\begin{bmatrix} 1 & 0  \\[0.05in] 0 & b_2 \end{bmatrix}^{-1}
\begin{bmatrix} \bar{\eta}^1 \\[0.05in] \bar{\eta}^2  \end{bmatrix} ,
\end{gather}
with $b_2  \neq 0$.

The 0-adapted coframings are the local sections of a principal f\/iber bundle
$\pi: \scB_0 \to \scX$, with structure group $G_0 \subset GL(2)$ consisting of all invertible matrices of the form in \eqref{2-1-g0-group}.  The right action of $G_0$ on sections $\sigma:\scX \to \scB_0$ is given by $\sigma \cdot g = g^{-1} \sigma$.  (This is the reason for the inverse occurring in \eqref{2-1-g0-group}.)

There exist canonical 1-forms $\eta^1$, $\eta^2$ on
$\scB_0$ with the {\em reproducing property} that for any coframing $(\bar{\eta}^1, \bar{\eta}^2)$ given by a local
section $\sigma: \scX \to \scB_0$,
\[ \sigma^{\ast} (\eta^i) = \bar{\eta}^i. \]
These are referred to as the {\em semi-basic} forms on $\scB_0$.
A standard argument shows that there also exists a (non-unique) 1-form $\beta_2$ (referred to as a {\em pseudo-connection form} or,
more succinctly, a {\em connection form}), linearly independent from
the semi-basic forms, and a function $T^1_{12}$ on $\scB_0$ (referred
to as a {\em torsion function}) such that
\begin{gather}\label{2-1-g0-structure-eqs}
\begin{bmatrix} d\eta^1 \\[0.05in] d\eta^2  \end{bmatrix} =
- \begin{bmatrix} 0 & 0 \\[0.05in] 0 & \beta_2  \end{bmatrix}
\begin{bmatrix} \eta^1 \\[0.05in] \eta^2  \end{bmatrix} +
\begin{bmatrix} T^1_{12} \eta^1 \& \eta^2 \\[0.05in] 0  \end{bmatrix}.
\end{gather}
(See \cite{Gardner89} for details about the method of equivalence.)
These are the {\em structure equations} of the $G_0$-structure
$\scB_0$.  The semi-basic forms and the connection form together form a
local coframing on $\scB_0$.

We proceed with the method of equivalence by examining how the
function $T^1_{12}$ varies if we change from one 0-adapted coframing
to another.  A straightforward computation shows that under a
transformation of the form \eqref{2-1-g0-group}, we have
\begin{gather}
 \tilde{T}^1_{12} = b_2 T^1_{12}. \label{2-1-torsion-trans}
\end{gather}
Thus $T^1_{12}$ is a {\em relative invariant}: if it vanishes for any coframing at a point $x \in \scX$, then it vanishes for every coframing at $x$.

At this point, we need to divide into cases based on whether or not $T^1_{12}$ vanishes.

{\bf Case 1.} Suppose that $T^1_{12}=0$.  Then the structure equations \eqref{2-1-g0-structure-eqs} contain no local invariants, so no further adaptations can be made.  The tableau consisting of all matrices of the form in~\eqref{2-1-g0-structure-eqs} is involutive with Cartan character $s_1=1$.

Given any 0-adapted coframing, consider the dual framing $(v_1, v_2)$.  The structure equations imply that
\[ [v_1, v_2] \equiv 0 \mod{v_2}; \]
therefore, $\scF$ is almost bracket-generating.  In order to recover the normal form of Theorem \ref{2-1-thm}, let $(\bar{\eta}^1, \bar{\eta}^2)$ be any 0-adapted local coframing on $\scX$ given by a local section $\sigma: \scX \to \scB_0$, and consider the pullbacks of the structure equations \eqref{2-1-g0-structure-eqs} via $\sigma$:
\begin{gather*}
d\bar{\eta}^1   = 0 , \qquad 
d\bar{\eta}^2   = -\bar{\beta}_2 \& \bar{\eta}^2. \notag
\end{gather*}
(Note that the pullback $\bar{\beta}_2$ of the connection form $\beta_2$ is a 1-form on $\scX$, so it is no longer linearly independent from $\bar{\eta}^1$, $\bar{\eta}^2$.)

Since $\bar{\eta}^1$ is an exact 1-form, we can choose a local coordinate $x^1$ on $\scX$ such that $\bar{\eta}^1 = dx^1$.  Then, since $\bar{\eta}^2$ is an integrable 1-form, we can choose a local coordinate $x^2$ on $\scX$, functionally independent from $x^1$, such that
\[ \bar{\eta}^2 = \lambda\, dx^2 \]
for some nonvanishing function $\lambda$ on $\scX$.  By replacing the given coframing with the (also 0-adapted) coframing $(\bar{\eta}^1, \lambda^{-1} \bar{\eta}^2)$, we can assume that $\lambda \equiv 1$.  Our 1-adapted coframing now has the form
\[ \bar{\eta}^1 = dx^1, \qquad \bar{\eta}^2 = dx^2. \]
The dual framing is
\[ v_1 = \frac{\partial}{\partial x^1}, \qquad v_2 = \frac{\partial}{\partial x^2}, \]
and
\[ \scF = \frac{\partial}{\partial x^1} + \text{span}\left(\frac{\partial}{\partial x^2}\right). \]

{\bf Case 2.} Suppose that $T^1_{12} \neq 0$.  \eqref{2-1-torsion-trans} implies that, given any 0-adapted coframing $(\bar{\eta}^1, \bar{\eta}^2)$, we can perform a transformation of the form \eqref{2-1-g0-group} to arrive at a new 0-adapted coframing for which $T^1_{12} = 1$.  Such a coframing will be called {\em $1$-adapted}.

Any two 1-adapted coframings on $\scX$ vary by a transformation of the form
\begin{gather*}
\begin{bmatrix} \tilde{\bar{\eta}}^1 \\[0.05in] \tilde{\bar{\eta}}^2 \end{bmatrix} =
\begin{bmatrix} 1 & 0  \\[0.05in] 0 & 1  \end{bmatrix}^{-1}
\begin{bmatrix} \bar{\eta}^1 \\[0.05in] \bar{\eta}^2 \end{bmatrix}.
\end{gather*}
In other words, the 1-adapted coframings are the local sections of an $(e)$-structure $\scB_1 \subset \scB_0$.  The projection $\pi: \scB_1 \to \scX$ is a dif\/feomorphism, and there is a {\em unique} 1-adapted coframing $(\bar{\eta}^1, \bar{\eta}^2)$ on $\scX$.
This coframing has structure equations
\begin{gather}
d\bar{\eta}^1   = \bar{\eta}^1 \& \bar{\eta}^2 ,\label{2-1-case2-structure} \\
d\bar{\eta}^2   = T^2_{12} \, \bar{\eta}^1 \& \bar{\eta}^2 ,\notag
\end{gather}
for some function $T^2_{12}$ on $\scX$.

Consider the dual framing $(v_1, v_2)$ of a 1-adapted coframing.  The structure equations \eqref{2-1-case2-structure} imply that
\[ [v_1, v_2] \equiv -v_1 \mod{v_2}; \]
therefore, $\scF$ is bracket-generating.  Since $\bar{\eta}^1$ is integrable but not exact, we can choose local coordinates $x^1$, $x^2$ on $\scX$ such that $\bar{\eta}^1 = \frac{1}{x^2} dx^1$.  Then the f\/irst equation in \eqref{2-1-case2-structure} implies that
\[ \bar{\eta}^2 = \frac{1}{x^2}\, (dx^2 - J\, dx^1) \]
for some function $J$ on $\scX$.  The second equation then implies that
\[ T^2_{12} = x^2\,J_{x^2} - J. \]
In terms of these local coordinates, the dual framing is
\[ v_1 = x^2 \left( \frac{\partial}{\partial x^1} + J \frac{\partial}{\partial x^2} \right), \qquad v_2 = x^2\, \frac{\partial}{\partial x^2}, \]
and since span$(v_2)$ = span$(\frac{\partial}{\partial x^2})$,
\begin{gather}
 \scF = x^2 \left( \frac{\partial}{\partial x^1} + J \frac{\partial}{\partial x^2} \right) + \text{span}\left(\frac{\partial}{\partial x^2}\right). \label{2-1-normal-form}
\end{gather}
\end{proof}

\begin{Remark}
This choice of coordinates $x^1$, $x^2$ such that $\bar{\eta}^1 = \frac{1}{x^2} dx^1$ is unique up to coordinate transformations of the form
\begin{gather}
 x^1 \to f(\tilde{x}^1), \qquad x^2 \to f'(\tilde{x}^1) \tilde{x}^2. \label{2-1-coord-change}
\end{gather}
Under a transformation of the form \eqref{2-1-coord-change}, we have
\[ J(x^1, x^2) \to J(f(\tilde{x}^1), f'(\tilde{x}^1) \tilde{x}^2) - \frac{f''(\tilde{x}^1)}{f'(\tilde{x}^1)} \tilde{x}^2. \]
Therefore, an admissible change of coordinates can only change $J$ by a linear function of $x^2$, and point-af\/f\/ine distributions of this type locally depend essentially on one arbitrary function of 2 variables.  In particular, the point-af\/f\/ine distribution~\eqref{2-1-normal-form} is locally equivalent to the ``f\/lat'' case $J=0$ if and only if $J = g(x^1) x^2$ for some function $g(x^1)$.
\end{Remark}

\subsection[Equivalence problem for rank $(n-1)$ point-affine distributions on $n$-manifolds]{Equivalence problem for rank $\boldsymbol{(n-1)}$ point-af\/f\/ine distributions\\ on $\boldsymbol{n}$-manifolds}

Theorem \ref{2-1-thm} is a special case of the following more general theorem:

\begin{Theorem}\label{n-n-1-thm}
Let $\scF$ be a rank $(n-1)$ strictly affine point-affine distribution of constant type on an $n$-dimensional manifold $\scX$.
Let $L_{\scF}$ be the associated direction distribution of rank $(n-1)$, let $v_1$ denote the distinguished vector field, and let $\bar{\eta}^1$ be the unique 1-form on $\scX$ satisfying
\[ L_{\scF} = (\bar{\eta}^1)^{\perp}, \qquad \bar{\eta}^1(v_1) = 1. \]
Let $k$ denote the Pfaff rank of $\bar{\eta}^1$, i.e., the unique integer such that
\[ \bar{\eta}^1 \& (d\bar{\eta}^1)^k \neq 0, \qquad \bar{\eta}^1 \& (d\bar{\eta}^1)^{k+1} = 0. \]
\begin{enumerate}\itemsep=0pt
\item If $(d\bar{\eta}^1)^{k+1} = 0$, then in a sufficiently small neighborhood of any point $x \in \scX$ on which $k$ is constant, there exist local coordinates $(x^1,\ldots, x^n)$ such that
\begin{gather*}
 \scF =   \left(1 + \sum_{r=1}^k x^{k+r+1} J_{k+r+1} \right) \frac{\partial}{\partial x^1} + \sum_{r=1}^k \left( J_{k+r+1} \frac{\partial}{\partial x^{r+1}} - J_{r+1} \frac{\partial}{\partial x^{k+r+1}} \right) \\
\phantom{\scF =}{}  + \text{\rm span}\left( \left(\frac{\partial}{\partial x^2} + x^{k+2} \frac{\partial}{\partial x^1} \right), \ldots, \left(\frac{\partial}{\partial x^{k+1}} + x^{2k+1} \frac{\partial}{\partial x^1}\right), \,\frac{\partial}{\partial x^{k+2}}, \ldots, \frac{\partial}{\partial x^n} \right),
\end{gather*}
where $J_2, \ldots, J_{2k+1}$ are arbitrary functions on $\scX$.
\item If $(d\bar{\eta}^1)^{k+1} \neq 0$, then in a sufficiently small neighborhood of any point $x \in \scX$ on which $k$ is constant, there exist local coordinates $(x^1,\ldots, x^n)$ such that
\begin{gather*}
 \scF =   \left(x^{2k+2}\! +\! \sum_{r=1}^k x^{k\!+\!r\!+\!1} J_{k\!+\!r\!+\!1} \right) \frac{\partial}{\partial x^1}\! +\! \sum_{r=1}^k \left( J_{k\!+\!r\!+\!1} \frac{\partial}{\partial x^{r+1}}\! - \! J_{r+1} \frac{\partial}{\partial x^{k\!+\!r\!+\!1}} \right)\! -\! J_1 \frac{\partial}{\partial x^{2k+2}} \\
\phantom{\scF = }{}
 + \text{\rm span}\left( \left(x^{2k+2} \left( \frac{\partial}{\partial x^2} + x^{k+2} \frac{\partial}{\partial x^1} \right) - J_2 \frac{\partial}{\partial x^{2k+2}} \right), \ldots,  \right. \\
\phantom{\scF = }{} \qquad \qquad
\left.  \left(x^{2k+2} \left( \frac{\partial}{\partial x^{k+1}} + x^{2k+1} \frac{\partial}{\partial x^1} \right) - J_{k+1} \frac{\partial}{\partial x^{2k+2}} \right), \right. \\
\phantom{\scF = }{} \qquad \qquad
\left.  \left(\frac{\partial}{\partial x^{k+2}} - J_{k+2} \frac{\partial}{\partial x^{2k+2}} \right), \ldots,  \left( \frac{\partial}{\partial x^{2k+1}} - J_{2k+1} \frac{\partial}{\partial x^{2k+2}} \right) ,  \right. \\
\phantom{\scF = }{} \qquad \qquad  \quad
\left.
  \frac{\partial}{\partial x^{2k+2}}, \ldots, \frac{\partial}{\partial x^n}
\right),
\end{gather*}
where $J_1, \ldots, J_{2k+1}$ are arbitrary functions on $\scX$.
\end{enumerate}
\end{Theorem}

\begin{proof}
Let $\scF$ be a rank $(n-1)$ point-af\/f\/ine distribution of constant type on an $n$-dimensional manifold $\scX$. A local framing $(v_1, \ldots, v_n)$ on $\scX$ will be called {\em admissible} if $v_1$ is the distinguished vector f\/ield and
\[ \scF = v_1 + \text{span}(v_2, \ldots, v_n). \]

A local coframing $(\bar{\eta}^1, \ldots, \bar{\eta}^n)$ on $\scX$ will be called {\em admissible} if, for any admissible framing $(v_1, \ldots, v_n)$,
\begin{gather*}
\bar{\eta}^1(v_1) = 1, \\ \bar{\eta}^1(v_j) = 0, \qquad 2 \leq j \leq n, \\
\bar{\eta}^j(v_1) = 0, \qquad 2 \leq j \leq n.
\end{gather*}
Any two admissible coframings on $\scX$ vary by a transformation of the form
\[
\begin{bmatrix} \tilde{\bar{\eta}}^1 \\ \tilde{\bar{\eta}}^2 \\ \vdots \\ \tilde{\bar{\eta}}^n \end{bmatrix} =
\left[
\begin{MAT}[2.5pt]{c|ccc}
 1 & 0 & \cdots & 0 \\- 0 &&&  \\ \vdots &&  A & \\ 0 &&&  \\
\end{MAT}
   \right]^{-1}
\begin{bmatrix} \bar{\eta}^1 \\ \bar{\eta}^2 \\ \vdots \\ \bar{\eta}^n \end{bmatrix}
\]
for some matrix $A \in GL(n-1, \R)$.

An admissible coframing will be called {\em $0$-adapted} if
\[ d\bar{\eta}^1 \equiv \sum_{r=1}^k \bar{\eta}^{r+1} \& \bar{\eta}^{k+r+1} \mod{\bar{\eta}^1}, \]
where $k$ is the Pfaf\/f rank of $\bar{\eta}^1$.
Since the Cartan system
\[ \calC(\bar{\eta}^1) = \{\bar{\eta}^1, \ldots, \bar{\eta}^{2k+1} \} \]
and the system
\[ v_1^{\perp} = \{\bar{\eta}^2, \ldots, \bar{\eta}^n \} \]
are both well-def\/ined independent of the choice of 0-adapted coframing, the subsystem
\[ \calC(\bar{\eta}^1) \cap v_1^{\perp} =  \{\bar{\eta}^2, \ldots, \bar{\eta}^{2k+1} \}  \]
is also well-def\/ined.  Thus any two 0-adapted coframings on $\scX$ vary by a transformation of the form
\begin{gather}
\begin{bmatrix}
\tilde{\bar{\eta}}^1 \\ \tilde{\bar{\eta}}^2 \\ \vdots \\ \tilde{\bar{\eta}}^{2k+1} \\ \tilde{\bar{\eta}}^{2k+2} \\ \vdots \\   \tilde{\bar{\eta}}^n \end{bmatrix}
\left[
\begin{MAT}[3pt]{c|ccc|ccc}
1 & 0 & \cdots & 0 & 0 & \cdots & 0 \\-
0 & & & &  0 & \cdots & 0 \\
\vdots & & A & & \vdots &  & \vdots  \\
0 & & & & 0 & \cdots & 0 \\-
0 & & & & & &   \\
\vdots & & B & & & C &  \\
0 & & & & & & \\
\end{MAT}
   \right]^{-1}
\begin{bmatrix} \bar{\eta}^1 \\ \bar{\eta}^2 \\ \vdots \\ \bar{\eta}^{2k+1} \\ \bar{\eta}^{2k+2} \\ \vdots \\   \bar{\eta}^n \end{bmatrix}, \label{n-n-1-g0-group}
\end{gather}
where $A \in Sp(2k, \R),\ B \in M_{(n-2k-1)\times(2k)}(\R),\ C \in GL(n\!-\!2k\!-\!1, \R)$.

The bundle $\pi: \scB_0 \to \scX$ of 0-adapted coframings with structure group $G_0$ consisting of matrices of the form \eqref{n-n-1-g0-group} is def\/ined as in Section~\ref{2-1-results}, as are the semi-basic forms $\eta^1, \ldots, \eta^n$ on~$\scB_0$.  The f\/irst structure equation of $\scB_0$ is
\begin{gather*}
d\eta^1 = \eta^1 \& \phi + \sum_{r=1}^k \eta^{r+1} \& \eta^{k+r+1}, 
\end{gather*}
where $\phi \in \pi^*(v_1^{\perp})$.  There are two cases to consider:
\begin{enumerate}\itemsep=0pt
\item $\phi \in \calC(\eta^1)$;
\item $\phi \not\in \calC(\eta^1)$.
\end{enumerate}

{\bf Case 1: $\boldsymbol{\phi \in \calC(\eta^1)}$.} 
Then
\[ (d\eta^1)^{k+1}=0. \]
Let $\sigma: \scX \to \scB_0$ be any local section of $\scB_0$ (i.e., a 0-adapted coframing $(\bar{\eta}^1, \ldots, \bar{\eta}^n)$ on $\scX$).  By the reproducing property of the semi-basic forms,
\[ (d\bar{\eta}^1)^{k+1}=0, \]
so by the Pfaf\/f theorem there exist local coordinates $x^1, \ldots, x^{2k+1}$ on $\scX$ such that
\begin{gather}
\bar{\eta}^1 = dx^1 - \sum_{r=1}^k x^{k+r+1}\, dx^{r+1}. \label{n-n-1-case1-eta1}
\end{gather}
We can write
\[ \bar{\phi} = \sum_{j=1}^{2k} J_{j+1} \bar{\eta}^{j+1} \]
for some functions $J_2, \ldots, J_{2k+1}$ on $\scX$.  Then we have
\begin{gather*}
d\bar{\eta}^1   = \bar{\eta}^1 \& \left( \sum_{j=1}^{2k} J_{j+1} \bar{\eta}^{j+1} \right) + \sum_{r=1}^k \bar{\eta}^{r+1} \& \bar{\eta}^{k+r+1} \\
\phantom{d\bar{\eta}^1}{}   = \sum_{r=1}^k \big( \bar{\eta}^{r+1} + J_{k+r+1} \bar{\eta}^1 \big) \& \big( \bar{\eta}^{k+r+1} - J_{r+1} \bar{\eta}^1 \big)   = \sum_{r=1}^k dx^{r+1} \& dx^{k+r+1}.
\end{gather*}
Thus we can choose a 0-adapted coframing with
\begin{gather}
\left.
\begin{array}{l}
 \bar{\eta}^{s+1}  = dx^{s+1} - J_{k+s+1} \bar{\eta}^1 \\
 \bar{\eta}^{k+s+1}  = dx^{k+s+1} + J_{s+1} \bar{\eta}^1
\end{array} \right\}, \ \ \ 1 \leq s \leq k. \label{n-n-1-case1-etas}
\end{gather}
For reasons of dimension, the system
\[ v_1^{\perp} = \{\bar{\eta}^2, \ldots, \bar{\eta}^n \} \]
is Frobenius; therefore, we can extend the local coordinates $x^1, \ldots, x^{2k+1}$ to a full local coordinate system $x^1, \ldots, x^n$ with the property that
\[ v_1^{\perp} = \{\bar{\eta}^2, \ldots, \bar{\eta}^{2k+1}, dx^{2k+2}, \ldots,  dx^n \}; \]
we can then choose the remainder of our 0-adapted coframing to be
\begin{gather}
 \bar{\eta}^j = dx^j, \qquad 2k+2 \leq j \leq n. \label{n-n-1-case1-etaj}
\end{gather}

The dual framing to the coframing def\/ined by \eqref{n-n-1-case1-eta1}, \eqref{n-n-1-case1-etas}, and \eqref{n-n-1-case1-etaj} is:
\begin{gather*}
v_1  = \left(1 + \sum_{r=1}^k x^{k+r+1} J_{k+r+1} \right) \frac{\partial}{\partial x^1} + \sum_{r=1}^k \left( J_{k+r+1} \frac{\partial}{\partial x^{r+1}} - J_{r+1} \frac{\partial}{\partial x^{k+r+1}} \right) ,\\
v_s  = \frac{\partial}{\partial x^s} + x^{k+s} \frac{\partial}{\partial x^1}, \qquad 2 \leq s \leq k+1, \\
v_j  = \frac{\partial}{\partial x^j}, \qquad k+2 \leq j \leq n.
\end{gather*}
The corresponding point-af\/f\/ine distribution $\scF$ is:
\begin{gather*}
 \scF =   \left(1 + \sum_{r=1}^k x^{k+r+1} J_{k+r+1} \right) \frac{\partial}{\partial x^1} + \sum_{r=1}^k \left( J_{k+r+1} \frac{\partial}{\partial x^{r+1}} - J_{r+1} \frac{\partial}{\partial x^{k+r+1}} \right) \\
\phantom{\scF =}{}  + \text{span}\left( \left(\frac{\partial}{\partial x^2} + x^{k+2} \frac{\partial}{\partial x^1} \right), \ldots, \left(\frac{\partial}{\partial x^{k+1}} + x^{2k+1} \frac{\partial}{\partial x^1}\right), \,\frac{\partial}{\partial x^{k+2}}, \ldots, \frac{\partial}{\partial x^n} \right).
\end{gather*}

Note that our choices of local coordinates were arbitrary up to at most functions of $(n-1)$ variables, while the local invariants $J_2, \ldots, J_{2k+1}$ are arbitrary functions of $n$ variables.  Therefore, point-af\/f\/ine distributions of this type locally depend essentially on $2k$ arbitrary functions of $n$ variables.

{\bf Case 2: $\boldsymbol{\phi \not\in \calC(\eta^1)}$.}
Then
\[ (d\eta^1)^{k+1} \neq 0. \]
Let $\sigma: \scX \to \scB_0$ be any local section of $\scB_0$ (i.e., a 0-adapted coframing $(\bar{\eta}^1, \ldots, \bar{\eta}^n)$ on $\scX$).  Because $d\eta^1$ is semi-basic, it follows from the reproducing property of the semi-basic forms that
\[ (d\bar{\eta}^1)^{k+1} \neq 0. \]
By the Pfaf\/f theorem, there exist local coordinates $x^1, \ldots, x^{2k+2}$ on $\scX$ such that
\begin{gather}
\bar{\eta}^1 = \frac{1}{x^{2k+2}} \left( dx^1 - \sum_{r=1}^k x^{k+r+1}\, dx^{r+1} \right). \label{n-n-1-case2-eta1}
\end{gather}
Then we have
\begin{gather*}
d\bar{\eta}^1   = \bar{\eta}^1 \& \bar{\phi} + \sum_{r=1}^k \bar{\eta}^{r+1} \& \bar{\eta}^{k+r+1} \\
\phantom{d\bar{\eta}^1 }{}  =   \bar{\eta}^1 \& \left( \frac{dx^{2k+2}}{x^{2k+2}} \right) + \frac{1}{x^{2k+2}} \sum_{r=1}^k dx^{r+1} \& dx^{k+r+1}.
\end{gather*}
Therefore,
\[ \sum_{r=1}^k \bar{\eta}^{r+1} \& \bar{\eta}^{k+r+1} \equiv \frac{1}{x^{2k+2}} \sum_{r=1}^k dx^{r+1} \& dx^{k+r+1} \mod{\bar{\eta}^1}. \]
Thus we can choose a 0-adapted coframing with
\begin{gather}
\left.
\begin{array}{l}
 \bar{\eta}^{s+1}   = \frac{1}{x^{2k+2}}(dx^{s+1} - J_{k+s+1} \bar{\eta}^1) \\
 \bar{\eta}^{k+s+1}   = dx^{k+s+1} + J_{s+1} \bar{\eta}^1
\end{array} \right\}, \ \ \ 1 \leq s \leq k, \label{n-n-1-case2-etas}
\end{gather}
for some functions $J_2, \ldots, J_{2k+1}$ on $\scX$.  Then we have
\[ \bar{\phi} \equiv \frac{1}{x^{2k+2}} \left( dx^{2k+2} + \sum_{j=1}^{2k} J_{j+1} \bar{\eta}^{j+1} \right) \mod{\bar{\eta}^1}. \]
The requirement that $\phi \in \pi^*(v_1^{\perp})$ is equivalent to $\bar{\phi} \in v_1^{\perp}$, and this determines a unique function~$J_1$ on~$\scX$ such that
\[ \bar{\phi} = \frac{1}{x^{2k+2}} \left( dx^{2k+2} + \sum_{j=0}^{2k} J_{j+1} \bar{\eta}^{j+1} \right). \]
We can now set
\begin{gather}
 \bar{\eta}^{2k+2} = \bar{\phi} = \frac{1}{x^{2k+2}} \left( dx^{2k+2} + \sum_{j=0}^{2k} J_{j+1} \bar{\eta}^{j+1} \right), \label{n-n-1-case2-eta2k2}
\end{gather}
and a similar argument to that given in the previous case shows that we can extend the local coordinates $x^1, \ldots, x^{2k+2}$ to a full local coordinate system $x^1, \ldots, x^n$ and choose the remainder of our 0-adapted coframing to be
\begin{gather}
 \bar{\eta}^j = dx^j, \qquad 2k+3 \leq j \leq n. \label{n-n-1-case2-etaj}
\end{gather}

The dual framing to the coframing def\/ined by \eqref{n-n-1-case2-eta1}, \eqref{n-n-1-case2-etas}, \eqref{n-n-1-case2-eta2k2}, and \eqref{n-n-1-case2-etaj} is:
\begin{gather*}
v_1  = \left(x^{2k+2}\! +\! \sum_{r=1}^k x^{k+r+1} J_{k+r+1} \right) \frac{\partial}{\partial x^1}\! +\! \sum_{r=1}^k \left( J_{k+r+1} \frac{\partial}{\partial x^{r+1}} \!- \! J_{r+1} \frac{\partial}{\partial x^{k+r+1}} \right)\! -\! J_1 \frac{\partial}{\partial x^{2k+2}} ,\\
v_s  = x^{2k+2} \left( \frac{\partial}{\partial x^s} + x^{k+s} \frac{\partial}{\partial x^1} \right) - J_s \frac{\partial}{\partial x^{2k+2}}, \qquad 2 \leq s \leq k+1, \\
v_s  = \frac{\partial}{\partial x^s} - J_s \frac{\partial}{\partial x^{2k+2}}, \qquad k+2 \leq s \leq 2k+1, \\
v_{2k+2}  = x^{2k+2} \frac{\partial}{\partial x^{2k+2}}, \qquad
v_j  = \frac{\partial}{\partial x^j}, \qquad 2k+3 \leq j \leq n.
\end{gather*}

Since span$(v_{2k+2})$ = span$(\frac{\partial}{\partial x^{2k+2}})$, the corresponding point-af\/f\/ine distribution $\scF$ is:
\begin{gather*}
 \scF =   \left(x^{2k+2}\! +\! \sum_{r=1}^k x^{k+r+1} J_{k+r+1} \right) \frac{\partial}{\partial x^1}\! +\! \sum_{r=1}^k \left( J_{k+r+1} \frac{\partial}{\partial x^{r+1}}\! -\! J_{r+1} \frac{\partial}{\partial x^{k+r+1}} \right)\! -\! J_1 \frac{\partial}{\partial x^{2k+2}} \\
\phantom{\scF =}{}  + \text{\rm span}\left( \left(x^{2k+2} \left( \frac{\partial}{\partial x^2} + x^{k+2} \frac{\partial}{\partial x^1} \right) - J_2 \frac{\partial}{\partial x^{2k+2}} \right), \ldots,  \right. \\
\phantom{\scF =}{}   \qquad \quad
\left.  \left(x^{2k+2} \left( \frac{\partial}{\partial x^{k+1}} + x^{2k+1} \frac{\partial}{\partial x^1} \right) - J_{k+1} \frac{\partial}{\partial x^{2k+2}} \right), \right. \\
\phantom{\scF =}{}  \qquad \quad
\left.  \left(\frac{\partial}{\partial x^{k+2}} - J_{k+2} \frac{\partial}{\partial x^{2k+2}} \right), \ldots,  \left( \frac{\partial}{\partial x^{2k+1}} - J_{2k+1} \frac{\partial}{\partial x^{2k+2}} \right) ,
  \frac{\partial}{\partial x^{2k+2}}, \ldots, \frac{\partial}{\partial x^n}
\right).
\end{gather*}

Note that our choices of local coordinates were arbitrary up to at most functions of $(n-1)$ variables, while the local invariants $J_1, \ldots, J_{2k+1}$ are arbitrary functions of $n$ variables.  Therefore, point-af\/f\/ine distributions of this type locally depend essentially on $2k+1$ arbitrary functions of $n$ variables.
\end{proof}

\begin{corollary}\label{3-2-cor}
Let $\scF$ be a rank $2$ strictly affine point-affine distribution of constant type on a~$3$-dimensional manifold~$\scX$.
\begin{enumerate}\itemsep=0pt
\item If $\scF$ is almost bracket-generating, then in a sufficiently small neighborhood of any point $x \in \scX$, there exist local coordinates $(x^1, x^2, x^3)$ such that
\[ \scF = \frac{\partial}{\partial x^1} + \text{\rm span} \left( \frac{\partial}{\partial x^2}, \frac{\partial}{\partial x^3} \right). \]
\item If $\scF$ is bracket-generating and $L_{\scF}$ is Frobenius, then in a sufficiently small neighborhood of any point $x \in \scX$, there exist local coordinates $(x^1, x^2, x^3)$ such that
\[ \scF = \left( x^2 \frac{\partial}{\partial x^1} - J_1 \frac{\partial}{\partial x^2}  \right)  + \text{\rm span} \left( \frac{\partial}{\partial x^2}, \frac{\partial}{\partial x^3} \right),  \]
where $J_1$ is an arbitrary function on $\scX$.
\item If $\scF$ is bracket-generating and $L_{\scF}$ is a contact distribution, then in a sufficiently small neighborhood of any point $x \in \scX$, there exist local coordinates $(x^1, x^2, x^3)$ such that
\[ \scF = \left( (1 + x^3 J_3) \frac{\partial}{\partial x^1} + J_3 \frac{\partial}{\partial x^2} - J_2 \frac{\partial}{\partial x^3} \right) + \text{\rm span} \left( x^3 \frac{\partial}{\partial x^1} + \frac{\partial}{\partial x^2}, \frac{\partial}{\partial x^3} \right), \]
where $J_2$, $J_3$ are arbitrary functions on $\scX$.
\end{enumerate}
\end{corollary}

\begin{proof}
Let $\bar{\eta}^1$ be as in Theorem \ref{n-n-1-thm}.

1. If $\scF$ is almost bracket-generating, then $d\bar{\eta}^1 = 0$, so $k=0$ and Case 1 of Theorem \ref{n-n-1-thm} applies.

2. If $\scF$ is bracket-generating and $L_{\scF}$ is Frobenius, then
\[ \bar{\eta}^1 \& d\bar{\eta}^1 = 0, \qquad d\bar{\eta}^1 \neq 0, \]
so $k=0$ and Case 2 of Theorem \ref{n-n-1-thm} applies.

3. If $\scF$ is bracket-generating and $L_{\scF}$ is a contact distribution, then
\[ \bar{\eta}^1 \& d\bar{\eta}^1 \neq 0, \qquad(d\bar{\eta}^1)^2 = 0, \]
so $k=1$ and Case 1 of Theorem \ref{n-n-1-thm} applies.
\end{proof}

\subsection{Equivalence problem for rank 1 point-af\/f\/ine distributions on 3-manifolds}

In this section we consider point-af\/f\/ine equivalence for rank 1 point-af\/f\/ine distributions $\scF$ on 3-dimensional manifolds.  We will assume that $\scF$ is bracket-generating or almost bracket-generating, since otherwise the equivalence problem may be reduced to a problem on a 2-dimensional manifold.

\begin{Theorem}\label{3-1-thm}
Let $\scF$ be a rank $1$ strictly affine point-affine distribution of constant type on a~$3$-dimensional manifold $\scX$.
\begin{enumerate}\itemsep=0pt
\item\label{3-1-thm-case1} If $\scF$ is almost bracket-generating, then in a sufficiently small neighborhood of any point $x \in \scX$, there exist local coordinates $(x^1, x^2, x^3)$ such that
\[ \scF = \left( \frac{\partial}{\partial x^1} + x^3 \frac{\partial}{\partial x^2} + J \frac{\partial}{\partial x^3} \right) + \text{\rm span} \left(\frac{\partial}{\partial x^3} \right), \]
where $J$ is an arbitrary function on  $\scX$.
\item \label{3-1-thm-case2}If $\scF$ is bracket-generating and $L_{\scF^2}$ $($which must have rank~$2)$ is Frobenius, then in a sufficiently small neighborhood of any point $x \in \scX$, there exist local coordinates $(x^1, x^2, x^3)$ such that
\[ \scF = \left( x^2 \frac{\partial}{\partial x^1} + x^3 \frac{\partial}{\partial x^2} + J \frac{\partial}{\partial x^3} \right) + \text{\rm span} \left(\frac{\partial}{\partial x^3} \right), \]
where $J$ is an arbitrary function on  $\scX$.
\item\label{3-1-thm-case3} If $\scF$ is bracket-generating and $L_{\scF^2}$ is not Frobenius, then in a sufficiently small neighborhood of any point $x \in \scX$, there exist local coordinates $(x^1, x^2, x^3)$ such that
\[ \scF = \left( \frac{\partial}{\partial x^1}  + J \left(x^3 \frac{\partial}{\partial x^1} + \frac{\partial}{\partial x^2} + H \frac{\partial}{\partial x^3}   \right) \right)  + \text{\rm span} \left( x^3  \frac{\partial}{\partial x^1} +  \frac{\partial}{\partial x^2}  + H  \frac{\partial}{\partial x^3} \right),  \]
where $H$, $J$ are arbitrary functions on $\scX$  satisfying $\frac{\partial H}{\partial x^1} \neq 0$.

\end{enumerate}
\end{Theorem}

Compare Theorem \ref{3-1-thm} with Theorem \ref{elkin-3d-thm}; in each case there is a new invariant $J$ which ref\/lects the restriction that $v_1$ is a f\/ixed vector f\/ield in $\scF$.

\begin{proof}
A local framing of $(v_1, v_2, v_3)$ on $\scX$ will be called {\em admissible} if $v_1$ is the distinguished vector f\/ield and
\[ \scF = v_1 + \text{span}(v_2). \]
If $(v_1, v_2, v_3)$ is an admissible framing on $\scX$, then any other admissible framing $(\tilde{v}_1, \tilde{v}_2, \tilde{v}_3)$ has the form{\samepage
\begin{gather*}
\tilde{v}_1   = v_1, \qquad
\tilde{v}_2   = b_2 v_2, \qquad 
\tilde{v}_3   = a_3 v_1 + b_3 v_2 + c_3 v_3, \notag
\end{gather*}
for some functions $a_3$, $b_2$, $b_3$, $c_3$ on $\scX$ which satisfy $b_2 c_3 \neq 0$.}

A local coframing $(\bar{\eta}^1, \bar{\eta}^2, \bar{\eta}^3)$ on $\scX$ will be called {\em $0$-adapted} if it is the dual coframing of some admissible framing.  Any two 0-adapted coframings on $\scX$ vary by a transformation of the form{\samepage
\begin{gather}\label{3-1-g0-group}
\begin{bmatrix} \tilde{\bar{\eta}}^1 \\[0.05in] \tilde{\bar{\eta}}^2 \\[0.05in] \tilde{\bar{\eta}}^3 \end{bmatrix} =
\begin{bmatrix} 1 & 0 & a_3 \\[0.05in] 0 & b_2 & b_3 \\[0.05in] 0 & 0 & c_3 \end{bmatrix}^{-1}
\begin{bmatrix} \bar{\eta}^1 \\[0.05in] \bar{\eta}^2 \\[0.05in] \bar{\eta}^3 \end{bmatrix},
\end{gather}
with $b_2 c_3 \neq 0$.}

The 0-adapted coframings are the local sections of a principal f\/iber bundle $\scB_0 \to \scX$, with structure group $G_0 \subset GL(3)$ consisting of all invertible matrices of the form in \eqref{3-1-g0-group}.  The structure equations on $\scB_0$ take the form
\begin{gather*}
\begin{bmatrix} d\eta^1 \\[0.05in] d\eta^2 \\[0.05in] d\eta^3 \end{bmatrix} =
- \begin{bmatrix} 0 & 0 & \alp_3 \\[0.05in] 0 & \beta_2 & \beta_3 \\[0.05in] 0 & 0 & \gamma_3 \end{bmatrix}
\begin{bmatrix} \eta^1 \\[0.05in] \eta^2 \\[0.05in] \eta^3 \end{bmatrix} +
\begin{bmatrix} T^1_{12} \eta^1 \& \eta^2 \\[0.05in] 0 \\[0.05in] T^3_{12} \eta^1 \& \eta^2 \end{bmatrix}.
\end{gather*}

Under a transformation of the form \eqref{3-1-g0-group}, the torsion functions $T^1_{12}$, $T^3_{12}$ transform as follows:
\begin{gather}
\tilde{T}^1_{12}   = b_2 T^1_{12} - \frac{a_3 b_2}{c_3} T^3_{12}, \qquad
\tilde{T}^3_{12}   = \frac{b_2}{c_3} T^3_{12}. \label{3-1-g0-torsion-trans}  
\end{gather}
The function $T^3_{12}$ is a relative invariant, and the assumption that $\scF$ is (almost) bracket-generating implies that $T^3_{12} \neq 0$.  So, given any 0-adapted coframing $(\bar{\eta}^1, \bar{\eta}^2, \bar{\eta}^3)$, \eqref{3-1-g0-torsion-trans} implies that we can perform a transformation of the form \eqref{3-1-g0-group} to arrive at a new 0-adapted coframing for which $T^3_{12}=1$, $T^1_{12} = 0.$  Such a coframing will be called {\em $1$-adapted}.

Any two 1-adapted coframings on $\scX$ vary by a transformation of the form
\begin{gather}\label{3-1-g1-group}
\begin{bmatrix} \tilde{\bar{\eta}}^1 \\[0.05in] \tilde{\bar{\eta}}^2 \\[0.05in] \tilde{\bar{\eta}}^3 \end{bmatrix} =
\begin{bmatrix} 1 & 0 & 0 \\[0.05in] 0 & b_2 & b_3 \\[0.05in] 0 & 0 & b_2 \end{bmatrix}^{-1}
\begin{bmatrix} \bar{\eta}^1 \\[0.05in] \bar{\eta}^2 \\[0.05in] \bar{\eta}^3 \end{bmatrix},
\end{gather}
with $b_2 \neq 0$.  The 1-adapted coframings are the local sections of a
principal f\/iber bundle $\scB_1 \subset \scB_0$, with structure
group $G_1$ consisting of all matrices of the form in \eqref{3-1-g1-group}.
When restricted to $\scB_1$, the connection forms $\alpha_3$,
$\gamma_3 - \beta_2$ become semi-basic, thereby
introducing new torsion terms into the structure equations of
$\scB_1$.  By adding multiples of the semi-basic forms to the
connection forms $\beta_2$, $\beta_3$ so as to absorb as much of the torsion as possible,
we can arrange that the structure equations of $\scB_1$ take the form
\begin{gather*}
\begin{bmatrix} d\eta^1 \\[0.05in] d\eta^2 \\[0.05in] d\eta^3 \end{bmatrix} =
- \begin{bmatrix} 0 & 0 & 0 \\[0.05in] 0 & \beta_2 & \beta_3 \\[0.05in] 0 & 0 & \beta_2 \end{bmatrix}
\begin{bmatrix} \eta^1 \\[0.05in] \eta^2 \\[0.05in] \eta^3 \end{bmatrix} +
\begin{bmatrix} (T^1_{13} \eta^1 + T^1_{23} \eta^2) \& \eta^3 \\[0.05in] 0 \\[0.05in] \eta^1 \& \eta^2 + T^3_{13} \eta^1 \& \eta^3 \end{bmatrix}.
\end{gather*}

Under a transformation of the form \eqref{3-1-g1-group}, we have
\begin{gather*}
\tilde{T}^1_{13}   = b_2 T^1_{13}, \qquad 
\tilde{T}^1_{23}   = b_2^2 T^1_{23},\qquad 
\tilde{T}^3_{13}   = T^3_{13} + \frac{b_3}{b_2}. \notag
\end{gather*}
The functions $T^1_{13}$, $T^1_{23}$ are relative invariants; however, given any 1-adapted coframing $(\bar{\eta}^1, \bar{\eta}^2, \bar{\eta}^3)$, we can perform a transformation of the form \eqref{3-1-g1-group} to arrive at a new 1-adapted coframing for which $T^3_{13} = 0$.  Such a coframing will be called {\em $2$-adapted}.

Any two 2-adapted coframings on $\scX$ vary by a transformation of the form
\begin{gather}\label{3-1-g2-group}
\begin{bmatrix} \tilde{\bar{\eta}}^1 \\[0.05in] \tilde{\bar{\eta}}^2 \\[0.05in] \tilde{\bar{\eta}}^3 \end{bmatrix} =
\begin{bmatrix} 1 & 0 & 0 \\[0.05in] 0 & b_2 & 0 \\[0.05in] 0 & 0 & b_2 \end{bmatrix}^{-1}
\begin{bmatrix} \bar{\eta}^1 \\[0.05in] \bar{\eta}^2 \\[0.05in] \bar{\eta}^3 \end{bmatrix},
\end{gather}
with $b_2 \neq 0$.  The 2-adapted coframings are the local sections of a
principal f\/iber bundle $\scB_2 \subset \scB_1$, with structure
group $G_2$ consisting of all matrices of the form in \eqref{3-1-g2-group}.
When restricted to~$\scB_2$, the connection form $\beta_3$ becomes semi-basic, thereby
introducing new torsion terms into the structure equations of~$\scB_2$.  By adding multiples of the semi-basic forms to the
connection form~$\beta_2$ so as to absorb as much of the torsion as possible,
we can arrange that the structure equations of $\scB_2$ take the form
\begin{gather}\label{3-1-g2-structure}
\begin{bmatrix} d\eta^1 \\[0.05in] d\eta^2 \\[0.05in] d\eta^3 \end{bmatrix} =
- \begin{bmatrix} 0 & 0 & 0 \\[0.05in] 0 & \beta_2 & 0 \\[0.05in] 0 & 0 & \beta_2 \end{bmatrix}
\begin{bmatrix} \eta^1 \\[0.05in] \eta^2 \\[0.05in] \eta^3 \end{bmatrix} +
\begin{bmatrix} (T^1_{13} \eta^1 + T^1_{23} \eta^2) \& \eta^3 \\[0.05in] T^2_{13} \eta^1 \& \eta^3 \\[0.05in] \eta^1 \& \eta^2 \end{bmatrix}.
\end{gather}

Under a transformation of the form \eqref{3-1-g2-group}, we have
\begin{gather}
\tilde{T}^1_{13}   = b_2 T^1_{13}, \qquad
\tilde{T}^1_{23}   = b_2^2 T^1_{23}, \qquad
\tilde{T}^2_{13}   = T^2_{13}.\label{3-1-g2-torsion-trans} 
\end{gather}
Thus the functions $T^1_{13}$, $T^1_{23}$ are relative invariants, while $T^2_{13}$ is a well-def\/ined function on $\scX$.
As in Section~\ref{2-1-results}, we will divide into cases based on the vanishing/nonvanishing of $T^1_{23}$ and $T^1_{13}$.

{\bf Case 1.} Suppose that $T^1_{23} = T^1_{13}=0$.  Since $T^2_{13}$ is well-def\/ined on $\scX$, no further adaptations can be made based on the structure equations \eqref{3-1-g2-structure}.  The tableau consisting of all matrices of the form in \eqref{3-1-g2-structure} is not involutive, so in principle the next step in the method of equivalence would be to prolong the structure equations.  However, equations \eqref{3-1-g2-structure} already contain enough information to produce a local normal form for such structures, so we will proceed in this direction.

Given any 1-adapted coframing, consider the dual framing $(v_1, v_2, v_3)$. The structure equations \eqref{3-1-g2-structure} imply that
\begin{gather*}
 [v_1, v_2]   \equiv -v_3 \mod{v_2}, \qquad
 [v_2, v_3]   \equiv 0 \mod{v_2, v_3}, \qquad
 [v_1, v_3]   \equiv 0 \mod{v_2, v_3}.
\end{gather*}
The f\/irst equation implies that $L_{\scF^2}$ is spanned by $\{v_2, v_3\}$; the second implies that $L_{\scF^2}$ is a~Frobenius distribution, and the third implies that $\scF$ is almost bracket-generating.  Thus, this case corresponds to Case~\ref{3-1-thm-case1} of Theorem~\ref{3-1-thm}.

In order to recover the normal form of Theorem \ref{3-1-thm}, let $(\bar{\eta}^1, \bar{\eta}^2, \bar{\eta}^3)$ be any 1-adapted local coframing on $\scX$ given by a local section $\sigma: \scX \to \scB_2$, and consider the pullbacks of the structure equations \eqref{3-1-g2-structure} via $\sigma$:
\begin{gather}
d\bar{\eta}^1   = 0, \notag \\
d\bar{\eta}^2   = -\bar{\beta}_2 \& \bar{\eta}^2 + \bar{T}^2_{13} \bar{\eta}^1 \& \bar{\eta}^3, \label{3-1-case1-structure} \\
d\bar{\eta}^3   = -\bar{\beta}_2 \& \bar{\eta}^3 + \bar{\eta}^1 \& \bar{\eta}^2. \notag
\end{gather}

Since $\bar{\eta}^1$ is an exact 1-form, we can choose a local coordinate $x^1$ on $\scX$ such that $\bar{\eta}^1 = dx^1$.  Then, since $\bar{\eta}^3$ is a contact form with
\[ d\bar{\eta}^3 \equiv 0 \mod{\{\bar{\eta}^3, dx^1\}} , \]
there exist local coordinates $x^2$, $x^3$, functionally independent from $x^1$, and a nonvanishing function $\lambda$ on $\scX$ such that
\[ \bar{\eta}^3 = \lambda(dx^2 - x^3\, dx^1). \]
By replacing the given coframing with the (also 1-adapted) coframing $(\bar{\eta}^1,\, \lambda^{-1} \bar{\eta}^2,\, \lambda^{-1} \bar{\eta}^3)$, we can assume that $\lambda \equiv 1$.  The third equation in \eqref{3-1-case1-structure} then implies that
\[ \bar{\eta}^2 = dx^3 + B dx^1 + C (dx^2 - x^3\, dx^1) \]
for some functions $B$, $C$ on $\scX$, and that
\[ \bar{\beta}_2 = C dx^1 + D (dx^2 - x^3\, dx^1) \]
for some additional function $D$ on $\scX$.  Finally, the second equation in \eqref{3-1-case1-structure} implies that
\[ C = \tfrac{1}{2} B_{x^3}, \qquad D = \tfrac{1}{2} B_{x^3 x^3}. \]

Our 2-adapted coframing now has the form
\begin{gather}
\bar{\eta}^1  = dx^1, \notag \\
\bar{\eta}^2  = dx^3 + B\, dx^1 + \tfrac{1}{2} B_{x^3} (dx^2 - x^3\, dx^1) \label{3-1-case1-normal-form}, \\
\bar{\eta}^3  = dx^2 - x^3\, dx^1. \notag
\end{gather}
The dual framing is
\[ v_1 = \frac{\partial}{\partial x^1} + x^3  \frac{\partial}{\partial x^2} - B \frac{\partial}{\partial x^3}, \qquad v_2 = \frac{\partial}{\partial x^3}, \qquad v_3 = \frac{\partial}{\partial x^2} - \tfrac{1}{2} B_{x^3} \frac{\partial}{\partial x^3}; \]
setting $J=-B$, we have:
\[ \scF = v_1 + \text{span}(v_2) = \left( \frac{\partial}{\partial x^1} + x^3 \frac{\partial}{\partial x^2} + J \frac{\partial}{\partial x^3} \right) + \text{span} \left(\frac{\partial}{\partial x^3} \right). \]

Comparing the structure equations \eqref{3-1-case1-structure} with those of the coframing \eqref{3-1-case1-normal-form}, we see that
\[ T^2_{13} = J_{x^2} - \tfrac{1}{2}(J_{x^1 x^3} + x^3 J_{x^2 x^3} + J J_{x^3 x^3}) + \tfrac{1}{4} (J_{x^3})^2. \]

{\bf Case 2.} Suppose that $T^1_{23} = 0$,  $T^1_{13} \neq 0$. \eqref{3-1-g2-torsion-trans} implies that, given any 2-adapted coframing $(\bar{\eta}^1, \bar{\eta}^2, \bar{\eta}^3)$, we can perform a transformation of the form \eqref{3-1-g2-group} to arrive at a new 2-adapted coframing for which $T^1_{13}=1$.  Such a coframing will be called {\em $3$-adapted}.

Any two 3-adapted coframings on $\scX$ vary by a transformation of the form
\begin{gather*}
\begin{bmatrix} \tilde{\bar{\eta}}^1 \\[0.05in] \tilde{\bar{\eta}}^2 \\[0.05in] \tilde{\bar{\eta}}^3 \end{bmatrix} =
\begin{bmatrix} 1 & 0 & 0 \\[0.05in] 0 & 1 & 0 \\[0.05in] 0 & 0 & 1 \end{bmatrix}^{-1}
\begin{bmatrix} \bar{\eta}^1 \\[0.05in] \bar{\eta}^2 \\[0.05in] \bar{\eta}^3 \end{bmatrix}.
\end{gather*}
In other words, the 3-adapted coframings are the local sections of an $(e)$-structure $\scB_3 \subset \scB_2$.  The projection $\pi: \scB_3 \to \scX$ is a dif\/feomorphism, and there is a {\em unique} 3-adapted coframing $(\bar{\eta}^1, \bar{\eta}^2, \bar{\eta}^3)$ on $\scX$.

When restricted to $\scB_3$, the connection form $\beta_2$ becomes semi-basic. The f\/irst structure equation in \eqref{3-1-g2-structure} now takes the form
\[ d\eta^1 = \eta^1 \& \eta^3. \]
Dif\/ferentiating this equation yields
\[ \beta_2 \& \eta^1 \& \eta^3 = 0; \]
therefore, $\beta_2$ is a linear combination of $\eta^1$ and $\eta^3$, and the structure equations of $\scB_3$ take the form
\begin{gather}
d\eta^1   = \eta^1 \& \eta^3, \notag \\
d\eta^2   = T^2_{12} \eta^1 \& \eta^2 + T^2_{13} \eta^1 \& \eta^3 + T^2_{23} \eta^2 \& \eta^3 ,  \label{3-1-case2-structure} \\
d\eta^3   = \eta^1 \& \eta^2 + T^2_{12} \eta^1 \& \eta^3 . \notag
\end{gather}

Given any 3-adapted coframing, consider the dual framing $(v_1, v_2, v_3)$. The structure equations \eqref{3-1-case2-structure} imply that
\begin{gather*}
 [v_1, v_2]   \equiv -v_3 \mod{v_2}, \qquad
 [v_2, v_3]   \equiv 0 \mod{v_2}, \qquad
 [v_1, v_3]   \equiv -v_1 \mod{v_2}.
\end{gather*}
The f\/irst equation implies that $L_{\scF^2}$ is spanned by $\{v_2, v_3\}$; the second implies that $L_{\scF^2}$ is a~Frobenius distribution, and the third implies that $\scF$ is bracket-generating.  Thus, this case corresponds to Case~\ref{3-1-thm-case2} of Theorem~\ref{3-1-thm}.

In order to recover the normal form of Theorem~\ref{3-1-thm}, let $(\bar{\eta}^1, \bar{\eta}^2, \bar{\eta}^3)$ be the 3-adapted local coframing on $\scX$ given by the unique local section $\sigma:\scX \to \scB_3$, and consider the pullbacks of the structure equations \eqref{3-1-case2-structure} via $\sigma$. Since $\bar{\eta}^1$ is an integrable 1-form (i.e., $d\bar{\eta}^1 \equiv 0 \mod{\bar{\eta}^1}$) but not exact, we can choose local coordinates $x^1$, $x^2$ on $\scX$ such that $\bar{\eta}^1 = \frac{1}{x^2}\, dx^1$.  Then the f\/irst equation in \eqref{3-1-case2-structure} implies that
\[ \bar{\eta}^3 \equiv \frac{1}{x^2}\, dx^2 \mod{dx^1}. \]
Now, since $d\bar{\eta}^3 \not\equiv 0 \mod{\bar{\eta}^3}$, there must exist another local coordinate $x^3$, functionally independent from $x^1$, $x^2$, such that
\[ \bar{\eta}^3 = \frac{1}{x^2} \left( dx^2- \frac{x^3}{x^2} dx^1 \right) . \]
Moreover, the third equation in \eqref{3-1-case2-structure} implies that
\[ \bar{\eta}^2 = \frac{1}{x^2} \left( (dx^3  + \frac{1}{x^2} B\, dx^1 + \frac{1}{x^2} C\,  \left( dx^2- \frac{x^3}{x^2} dx^1 \right) \right)  \]
for some functions $B$, $C$ on $\scX$.  Finally, the second equation in \eqref{3-1-case2-structure} implies that
\[ C = \tfrac{1}{2}\left( x^2 B_{x^3} - x^3 \right). \]

Our 3-adapted coframing now has the form
\begin{gather}
\bar{\eta}^1    = \frac{1}{x^2} dx^1, \notag   \\
\bar{\eta}^2    = \frac{1}{x^2} dx^3 + \frac{1}{(x^2)^2} B\,dx^1 + \tfrac{1}{2}\left( \frac{(x^2 B_{x^3} - x^3)}{(x^2)^2} \right) \left( dx^2- \frac{x^3}{x^2} dx^1 \right) ,\label{3-1-case2-normal-form}  \\
\bar{\eta}^3    = \frac{1}{x^2}  dx^2 - \frac{x^3}{(x^2)^2} dx^1 . \notag
\end{gather}
The dual framing is
\begin{gather*}
v_1   = x^2 \frac{\partial}{\partial x^1} + x^3  \frac{\partial}{\partial x^2} - B \frac{\partial}{\partial x^3} ,\\
v_2   = x^2 \frac{\partial}{\partial x^3} ,\\
v_3   = x^2 \frac{\partial}{\partial x^2} - \tfrac{1}{2}\left( x^2 B_{x^3} - x^3 \right) \frac{\partial}{\partial x^3};
\end{gather*}
setting $J=-B$ and noting that $\text{span}(v_2) = \text{span}(\frac{\partial}{\partial x^3})$, we have:
\[ \scF = v_1 + \text{span}(v_2) = \left( x^2 \frac{\partial}{\partial x^1} + x^3 \frac{\partial}{\partial x^2} + J \frac{\partial}{\partial x^3} \right) + \text{span} \left(\frac{\partial}{\partial x^3} \right). \]

Comparing the structure equations \eqref{3-1-case2-structure} with those of the coframing \eqref{3-1-case2-normal-form}, we see that
\begin{gather*}
T^2_{12}   = \frac{1}{2x^2} (x^2 J_{x^3} - 3x^3) ,  \\
T^2_{13}   = \frac{1}{4 (x^2)^2} \Big{(} 3(x^3)^2  -  6 x^2 J  + 4 (x^2)^2 J_{x^2}  +  2 x^2 x^3 J_{x^3}  -  2 (x^2)^3 J_{x^1 x^3} \\
\phantom{T^2_{13}   =}{}  -  2 (x^2)^2 x^3 J_{x^2 x^3}  +  (x^2)^2 (J_{x^3})^2  -  2 (x^2)^2 J J_{x^3 x^3} \Big{)}, \\
T^2_{23}   =  \tfrac{1}{2}(1 - x^2 J_{x^3 x^3})   .
\end{gather*}

{\bf Case 3.} Suppose that $T^1_{23} \neq 0$. \eqref{3-1-g2-torsion-trans} implies that, given any 2-adapted coframing $(\bar{\eta}^1, \bar{\eta}^2, \bar{\eta}^3)$, we can perform a transformation of the form \eqref{3-1-g2-group} to arrive at a new 1-adapted coframing for which $T^1_{23}=\pm 1$.  Such a coframing will be called {\em $3$-adapted}.  (For ease of notation, let $\eps = T^1_{23} = \pm 1$.)

Any two 3-adapted coframings on $\scX$ vary by a transformation of the form
\begin{gather*}
\begin{bmatrix} \tilde{\bar{\eta}}^1 \\[0.05in] \tilde{\bar{\eta}}^2 \\[0.05in] \tilde{\bar{\eta}}^3 \end{bmatrix} =
\begin{bmatrix} 1 & 0 & 0 \\[0.05in] 0 & \pm 1 & 0 \\[0.05in] 0 & 0 & \pm 1 \end{bmatrix}^{-1}
\begin{bmatrix} \bar{\eta}^1 \\[0.05in] \bar{\eta}^2 \\[0.05in] \bar{\eta}^3 \end{bmatrix},
\end{gather*}
where the signs are chosen such that the matrix above has determinant equal to 1.  The corresponding principal f\/iber bundle $\scB_3 \subset \scB_2$ has discrete f\/iber group $G_3 \cong \mathbb{Z}/2\mathbb{Z}$, and we can regard~$\scB_3$ as an $(e)$-structure over a 2-to-1 cover of $\scX$.
When restricted to~$\scB_3$, the connection form $\beta_2$ becomes semi-basic; the structure equations of $\scB_3$ take the form
\begin{gather}
d\eta^1   = T^1_{13} \eta^1 \& \eta^3 + \eps \eta^2 \& \eta^3, \notag \\
d\eta^2   = T^2_{12} \eta^1 \& \eta^2 + T^2_{13} \eta^1 \& \eta^3 + T^2_{23} \eta^2 \& \eta^3 ,  \label{3-1-case3-structure} \\
d\eta^3   = \eta^1 \& \eta^2 + T^2_{12} \eta^1 \& \eta^3 + T^3_{23} \eta^2 \& \eta^3. \notag
\end{gather}

Given any 3-adapted coframing, consider the dual framing $(v_1, v_2, v_3)$. The structure equations \eqref{3-1-case3-structure} imply that
\begin{gather*}
 [v_1, v_2]   \equiv -v_3 \mod{v_2}, \qquad
 [v_2, v_3]   \equiv -\eps v_1 \mod{v_2, v_3}, \qquad
 [v_1, v_3]   \equiv 0 \mod{v_1, v_2, v_3}.
\end{gather*}
The f\/irst equation implies that $L_{\scF^2}$ is spanned by $\{v_2, v_3\}$, and the second implies that $L_{\scF^2}$ is not Frobenius and that $\scF$ is bracket-generating.  Thus, this case corresponds to Case \ref{3-1-thm-case3} of Theorem~\ref{3-1-thm}.

In order to recover the normal form of Theorem~\ref{3-1-thm}, let $(\bar{\eta}^1, \bar{\eta}^2, \bar{\eta}^3)$ be any 3-adapted local coframing on $\scX$ given by a local section $\sigma:\scX \to \scB_3$, and consider the pullbacks of the structure equations \eqref{3-1-case3-structure} via $\sigma$.
Since $d\bar{\eta}^1 \not\equiv 0 \mod{\bar{\eta}^1}$, Pfaf\/f's theorem implies that we can choose local coordinates $x^1$, $x^2$, $x^3$ on $\scX$ such that
\[ \bar{\eta}^1 = dx^1 - x^3\, dx^2. \]
Set $T^1_{13} = B$; this will simplify notation in what follows.
The f\/irst equation in \eqref{3-1-case3-structure} becomes
\[ ( B \bar{\eta}^1 + \eps \bar{\eta}^2) \& \bar{\eta}^3 = dx^2 \& dx^3, \]
which implies that
\begin{gather*}
\bar{\eta}^2  = \eps(-B (dx^1 - x^3\, dx^2) - \lambda^{-1}dx^2 + C (H dx^2 - dx^3)), \\
\bar{\eta}^3  = \lambda(H dx^2 - dx^3),
\end{gather*}
for some functions $\lambda$, $C$, $H$ on $\scX$ with $\lambda \neq 0$.  (If necessary, perform the contact transformation
\[ x^1 \to x^1 - x^2 x^3, \qquad x^2 \to x^3, \qquad x^3 \to -x^2 \]
to achieve this form.)
Now the third equation in \eqref{3-1-case3-structure} implies that
\[ -\eps \lambda^2 H_{x^1} = 1 . \]
Therefore, $H_{x^1} \neq 0$; $\eps H_{x^1} < 0$, and
\[ \lambda = \pm \frac{1}{\sqrt{-\eps H_{x^1}}}. \]
The $\mathbb{Z}/2\mathbb{Z}$ freedom in the structure group $G_3$ may be used to eliminate the sign ambiguity in $\lambda$; for simplicity we will take the positive square root.

Our 3-adapted coframing now has the form
\begin{gather}
\bar{\eta}^1    = dx^1 - x^3\, dx^2, \notag \\
\bar{\eta}^2    = -\eps \left(\sqrt{-\eps H_{x^1}} dx^2 + B (dx^1 - x^3\, dx^2) + C (H\, dx^2 - dx^3)\right)  ,
\label{3-1-case3-normal-form} \\
\bar{\eta}^3    = \frac{1}{\sqrt{-\eps H_{x^1}}} (H\, dx^2 - dx^3) . \notag
\end{gather}
The dual framing is
\begin{gather*}
v_1   =   \frac{\partial}{\partial x^1} - \frac{B}{\sqrt{-\eps H_{x^1}}} \left(x^3 \frac{\partial}{\partial x^1} + \frac{\partial}{\partial x^2} + H \frac{\partial}{\partial x^3}   \right),    \\
v_2   = -\frac{\eps}{\sqrt{-\eps H_{x^1}}} \left( x^3  \frac{\partial}{\partial x^1} +  \frac{\partial}{\partial x^2}  + H  \frac{\partial}{\partial x^3} \right),   \\
v_3   = -\sqrt{-\eps H_{x^1}}  \frac{\partial}{\partial x^3} + C \left(x^3 \frac{\partial}{\partial x^1} + \frac{\partial}{\partial x^2} + H \frac{\partial}{\partial x^3}   \right) ;
\end{gather*}
setting $J=  \frac{-B}{\sqrt{-\eps H_{x^1}}}$, we have:
\begin{gather*} \scF   = v_1 + \text{span}(v_2) \\
\phantom{\scF}{}  = \left( \frac{\partial}{\partial x^1}  + J \left(x^3 \frac{\partial}{\partial x^1} + \frac{\partial}{\partial x^2} + H \frac{\partial}{\partial x^3}   \right) \right)  + \text{span} \left( x^3  \frac{\partial}{\partial x^1} +  \frac{\partial}{\partial x^2}  + H  \frac{\partial}{\partial x^3} \right).
\end{gather*}

Comparing the structure equations \eqref{3-1-case3-structure} with those of the coframing \eqref{3-1-case3-normal-form} shows that
\begin{gather*}
 C = \frac{1}{2 H_{x^1}}\sqrt{-\eps H_{x^1}} \Bigg{(}(H_{x^1} H_{x^3} - H H_{x^1 x^3} - H_{x^1 x^2}) J \\
 \phantom{ C = \frac{1}{2 H_{x^1}}\sqrt{-\eps H_{x^1}} \Bigg{(}}{}
 - (1+x^3 J)H_{x^1 x^1} - (J_{x^2} + x^3 J_{x^1} + H J_{x^3}) H_{x^1}
  \Bigg{)},
\end{gather*}
and that
\begin{gather*}
T^2_{12}   =\tfrac{1}{2} \left(J_{x^2} + x^3 J_{x^1} + H J_{x^3} + J H_{x^3}\right) ,\\
T^2_{13}   = \eps C^2 H_{x^1} - H_{x^1} J_{x^3} + J^2 H_{x^1}  +  \frac{1}{2 \sqrt{-\eps H_{x^1}}} \bigg{(} 2C(H_{x^1} J_{x^2} + J H_{x^1, x^2} + H_{x^1 x^1} \\
\phantom{T^2_{13}   =}{}  + 2x^3 H_{x^1} J_{x^1} + x^3 J H_{x^1 x^1} + 2 H H_{x^1} J_{x^3} + H H_{x^1 x^3} J - 2 H_{x^1} H_{x^3} J) \\
\phantom{T^2_{13}   =}{}   - 2 x^3 C_{x^1} H_{x^1} J - 2 C_{x^1} H_{x^1} - 2 C_{x^3} H H_{x^1} J - 2 C_{x^2} H_{x^1} J
 \bigg{)}, \\
T^2_{23}   = -\left(C_{x^2} + x^3 C_{x^1} + H C_{x^3} + C H_{x^3}\right) - \frac{\eps}{2 \sqrt{-\eps H_{x^1}}} \left(H_{x^1 x^3} + 2J H_{x^1}\right), \\
T^3_{23}   = \frac{\eps}{2 H_{x^1} \sqrt{-\eps H_{x^1}}} \left(x^3 H_{x^1 x^1} - 2 H_{x^1} H_{x^3} + H_{x^1 x^2} + H H_{x^1 x^3}\right).\tag*{\qed}
\end{gather*}
\renewcommand{\qed}{}
\end{proof}

\section{Examples and questions for further study}\label{future}

\subsection{Examples}

\begin{Example}[NMR control]
In quantum mechanics, a particle with spin $\tfrac{1}{2}$ is represented by a state vector
\begin{gather*}
  | \psi \rangle = \begin{bmatrix}
                                   \alpha\\
                                   \beta\\
                                \end{bmatrix} \in \C^{2}
\end{gather*}
with $\ds{|\alpha|^{2} + |\beta|^{2} = 1}$.  In the presence of a magnetic f\/ield $\mbf{B} = B_{x}\mbf{i} + B_{y}\mbf{j} + B_{z}\mbf{k}$, the state vector at time $t$ is
\begin{gather*}
  | \psi(t) \rangle = U(t) | \psi_{0} \rangle,
\end{gather*}
where $U: \R \to SU(2)$ is a solution of Schr\"{o}dinger's equation
\begin{gather}\label{Schrodinger}
  \dot{U} = -i \calH\, U.
\end{gather}
(See \cite{NC00} for more details.) The Hamiltonian operator $\calH$ is given by
\begin{gather*}
  \calH = -\gamma \hbar (B_{x} \cdot I_{x} + B_{y} \cdot I_{y} + B_{z} \cdot I_{z}),
\end{gather*}
where $\gamma$ is the \emph{gyromagnetic ratio} (a constant), $\hbar = \frac{h}{2\pi}$ (where $h$ is Planck's constant), and~$I_{x}$,~$I_{y} $, and $I_{z}$ are the \emph{Pauli spin matrices}
\begin{gather*}
  I_{x} = \frac{1}{2} \begin{bmatrix}
                                    0 & 1\\
                                    1 & 0
                                  \end{bmatrix}, \qquad
 I_{y} = \frac{1}{2} \begin{bmatrix}
                                                                                             0 & -i\\
                                                                                             i & 0
                                                                                            \end{bmatrix}, \qquad I_{z} = \frac{1}{2} \begin{bmatrix}
                                                                                                                                                            1 & 0\\
                                                                                                                                                            0 & -1
                                                                                                                                                      \end{bmatrix}.
\end{gather*}
By multiplying the Pauli spin matrices by $-i$ as in~\eqref{Schrodinger} we obtain members of the Lie algebra of $SU(2)$.  The commutation relations for $-i I_{x}$, $-i I_{y}$, and $-i I_{z}$ are
\begin{gather*}
  \left[-i I_{x}, -i I_{y}\right]  = - i I_{z}, \qquad 
  \left[-i I_{y}, -i I_{z}\right]  = - i I_{x},\qquad 
  \left[-i I_{z}, -i I_{x}\right]  =  -i I_{y}.
\end{gather*}
In nuclear magnetic resonance (NMR) experiments, an ensemble of atomic nuclei of the same type can be modeled as a single spin system~\cite{EBW87}.  In such experiments, an ensemble of nuclei is subjected to a constant magnetic f\/ield $\mbf{B_{0}}$ in the $z$-direction (called the \emph{longitudinal field}), and an oscillating f\/ield $\mbf{B_{1}}$ in the $xy$-plane (called the \emph{rf field}):
\begin{gather*}
  \mbf{B_{0}}  = B_{0} \mbf{k},\\ 
  \mbf{B_{1}}  = B_{1}(\cos{(\omega t + \phi)} \mbf{i} + \sin{(\omega t + \phi)} \mbf{j}),
\end{gather*}
where $\omega$ is the frequency of the rf f\/ield and $\phi$ is a phase.  The Hamiltonian $\calH$ can thus be written~as
\begin{gather*}
  \calH = -\gamma \hbar (B_{0} \cdot I_{z} + B_{1}(\cos{(\omega t + \phi)} I_{x} + \sin{(\omega t + \phi)} I_{y}).
\end{gather*}
In NMR experiments, the magnitude of the longitudinal f\/ield $\mbf{B_{0}}$ is usually kept constant, while the magnitude of the rf f\/ield $\mbf{B_{1}}$ may be used as a control variable.  By transforming the coordinate system to a system that rotates at the same frequency $\omega$ as the rf f\/ield, and normalizing the constant $-\gamma \hbar B_{0}$ to 1, we may rewrite the Hamiltonian as
\begin{gather*}
  \calH = I_{z} + u(t) I_{x},
\end{gather*}
where $u(t) = \frac{B_{1}(t)}{B_{0}}$.

Thus Schr\"{o}dinger's equation \eqref{Schrodinger} becomes the control-af\/f\/ine system
\begin{gather}
 \dot{U} = -i I_z U - (i I_x U) u \label{example-4-1-system}
\end{gather}
on the 3-dimensional state space $SU(2)$.  This system corresponds to a strictly af\/f\/ine, rank 1 point-af\/f\/ine distribution
\[ \scF = v_1 + \text{span}(v_2), \]
where $v_1$, $v_2$ are the right-invariant vector f\/ields
\[ v_1 = -i I_z U, \qquad v_2 = -i I_x U \]
on $SU(2)$.  Let $\alpha_1$, $\alpha_2$, $\alpha_3$ denote the canonical right-invariant 1-forms on $SU(2)$ dual to the basis $-iI_x U$, $-iI_y U$, $-iI_z U$ for the right-invariant vector f\/ields; these forms satisfy the structure equations
\begin{gather*}
d\alpha_1  = \alpha_2 \& \alpha_3, \qquad
d\alpha_2  = \alpha_3 \& \alpha_1, \qquad
d\alpha_3  = \alpha_1 \& \alpha_2.
\end{gather*}
A 0-adapted coframing for the system \eqref{example-4-1-system} (corresponding to the choice $v_3 = -[v_1, v_2] = -iI_yU$) is given by
\[ \eta^1 = \alpha_3, \qquad \eta^2 = \alpha_1, \qquad \eta^3 = \alpha_2. \]
This coframing has structure equations
\begin{gather*}
d\eta^1 = \eta^2 \& \eta^3, \qquad
d\eta^2 = \eta^3 \& \eta^1, \qquad
d\eta^3 = \eta^1 \& \eta^2;
\end{gather*}
therefore this example falls into Case 3 of Theorem \ref{3-1-thm}, with $\eps=1$.  Since $T^1_{13}=0$, we have $J=0$.

Let $x^1$, $x^2$, $x^3$ be local coordinates such that $\eta^1 = dx^1 - x^3\, dx^2.$ Then as in the proof of Theorem \ref{3-1-thm}, there exist functions $C$, $H$, $\lambda$  such that
\begin{gather*}
\bar{\eta}^2  = -\lambda^{-1}dx^2 + C (H dx^2 - dx^3)), \\
\bar{\eta}^3  = \lambda(H dx^2 - dx^3).
\end{gather*}
(Note that $B=0$ in this example.)  The structure equations imply that
\begin{gather*}
 H = -e^{2f(x^2, x^3)} \tan\left(x^1 + g(x^2, x^3)\right) + h(x^2, x^3), \\
 C = -e^{-f(x^2,x^3)} \sin\left(x^1 + g(x^2, x^3)\right), \\
 \lambda = e^{-f(x^2,x^3)} \cos\left(x^1 + g(x^2, x^3)\right),
\end{gather*}
for functions $f$, $g$, $h$ which satisfy the additional PDEs
\begin{gather*}
e^{2f}f_{x^3} + g_{x^2} + h g_{x^3} + x^3 = 0, \\
e^{2f}g_{x^3} + f_{x^2} + h f_{x^3} - h_{x^3} = 0.
\end{gather*}
Thus the invariants for this system are
\[ H = -e^{2f(x^2, x^3)} \tan\left(x^1 + g(x^2, x^3)\right) + h(x^2, x^3), \qquad J=0. \]
\end{Example}

\begin{Example}[Navigation on a river]\label{boat-ex}
  Consider a boat navigating on a river with a current; for simplicity, assume that the current runs parallel to the $x$-axis with constant speed $c$.  The state of the boat is represented by its position $(x,y)$ and heading angle $\psi$ (see Fig.~\ref{boat-fig}); thus the state space is $\scX = \mathbb{R}^2 \times S^1$.

  {\setlength{\unitlength}{2pt}
\begin{figure}[th!]
\centering
\begin{picture}(40,40)(-20,-17)
\put(-15,0){\line(1,1){15}}
\put(0,-15){\line(1,1){15}}
\qbezier(0,15)(9,20)(17,17)
\qbezier(15,0)(20,9)(17,17)
\put(-15,0){\line(1,-3){3.6}}
\put(0,-15){\line(-3,1){11}}
\put(-20,-20){\line(1,1){45}}
\put(-30,3){\line(1,0){65}}
\qbezier(6,6)(7.3,5.6)(8,3)
\put(11,6){\makebox(0,0){$\psi$}}
\end{picture}
\caption{Heading angle $\psi$ in Example \ref{boat-ex}.}
\label{boat-fig}
\end{figure}
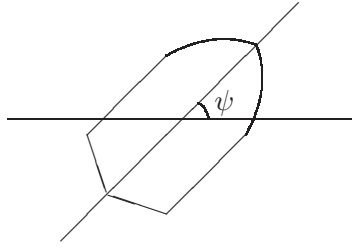}

Due to asymmetry in the shape of the boat's hull, the current may have a rotational ef\/fect as well as a translational ef\/fect on the boat, and this ef\/fect depends on the boat's heading angle.  Thus we will assume that the drift vector f\/ield has the form
\[ v_1 = c \frac{\partial}{\partial x} + r(\psi) \frac{\partial}{\partial \psi}. \]
We also assume that the boat can be propelled in a forward or backward direction, corresponding to the vector f\/ield
\[ v_2 = (\cos \psi) \frac{\partial}{\partial x} + (\sin \psi) \frac{\partial}{\partial y}, \]
and that it can be steered to the right or left, corresponding to the vector f\/ield
\[ v_3 = \frac{\partial}{\partial \psi}. \]
The navigation problem for the boat is then the control-af\/f\/ine system corresponding to the rank~2 point-af\/f\/ine distribution
\[ \scF = v_1 + \text{span}(v_2, v_3) \]
on $\scX$.

Note that $\scF$ fails to be strictly af\/f\/ine precisely when $\psi=0$ or $\psi=\pi$, since then $v_1$ is contained in $L_{\scF}$.  Thus our results are only applicable when the boat is {\em not} pointed directly upstream or downstream.  Since $L_{\scF} = \text{span}(v_2, v_3)$ is a contact distribution, this system falls into Case 3 of Corollary \ref{3-2-cor}.

The dual coframing to the framing $(v_1, v_2, v_3)$ on $\scX$ is
\[ \bar{\eta}^1 = \frac{1}{c} \left(dx - (\cot \psi) dy\right), \qquad \bar{\eta}^2 = (\csc \psi) dy, \qquad \bar{\eta}^3 = d\psi - \frac{r(\psi)}{c}\left(dx - (\cot \psi) dy\right). \]
This coframing has
\[ d\bar{\eta}^1 \equiv -\frac{(\csc \psi)}{c} \bar{\eta}^2 \& \bar{\eta}^3 \mod{\bar{\eta}^1}; \]
in order to obtain a 0-adapted coframing as in the proof of Theorem \ref{n-n-1-thm}, we multiply $\bar{\eta}^2$ by $\frac{(\sin \psi)}{c}$ and  $\bar{\eta}^3$ by $-(\csc^2 \psi)$, resulting in the following coframing:
\begin{gather}
 \bar{\eta}^1 = \frac{1}{c} \left(dx - (\cot \psi) dy\right), \qquad \bar{\eta}^2 = \frac{1}{c} dy, \nonumber\\ \bar{\eta}^3 = -(\csc^2 \psi) \left( d\psi - \frac{r(\psi)}{c}(dx - (\cot \psi) dy) \right). \label{ex-2-coframing-1}
\end{gather}
The local coordinate functions
\[ x^1 = \frac{x}{c}, \qquad x^2 = \frac{y}{c}, \qquad x^3 = \cot \psi \]
are Pfaf\/f normal form coordinates for $\bar{\eta}^1$; i.e.,
\[ \bar{\eta}^1 = dx^1 - x^3\, dx^2. \]
In terms of these coordinates, the coframing \eqref{ex-2-coframing-1} becomes:
\begin{gather*}
 \bar{\eta}^1 = dx^1 - x^3\, dx^2, \qquad \bar{\eta}^2 = dx^2, \qquad \bar{\eta}^3 = dx^3 -  ((x^3)^2+1) R(x^3)(dx^1 - x^3\, dx^2), 
\end{gather*}
where $R(x^3) = r(\psi) = r(\cot^{-1} (x^3))$.  Thus the invariants for this system are
\[ J_2 = ((x^3)^2+1) R(x^3) = (\csc^2 \psi) r(\psi), \qquad J_3=0. \]
\end{Example}

\subsection{Questions for further study}

There are two main issues that we hope to address in future papers, both motivated by optimal control theory.
\begin{enumerate}\itemsep=0pt

\item What metric structures (analogous to sub-Riemannian or sub-Finsler geometry for linear distributions) are appropriate for point-af\/f\/ine distributions, and what can we say about their geometry, geodesics, etc.?  Even for linear distributions, issues such as controllability and the presence of {\em rigid curves}~-- i.e., curves with no $C^1$ variations whatsoever among horizonal curves~-- are nontrivial, and the study of optimal trajectories for sub-Riemannian and sub-Finsler metrics is quite complicated.  (See, for instance, \cite{LS95, Montgomery02,CMW07}.)  Agrachev and Sarychev have studied time-optimal extremals and given necessary and suf\/f\/icient conditions for ridigity of trajectories for af\/f\/ine distributions in \cite{AS96}, but the consequences of imposing any type of metric structure on a~point-af\/f\/ine distribution and the resulting geometry remain entirely unexplored.  These issues are obviously important for understanding the associated control problems.

\item What can we say about the geometry of af\/f\/ine distributions of non-constant type, particularly for af\/f\/ine distributions which are not strictly af\/f\/ine?  Normal forms for non-strictly af\/f\/ine distributions of rank 1 have been given by Kang and Krener \cite{KK92}, Kang \cite{Kang96}, and Tall and Respondek \cite{TR02}, and for corank-1 af\/f\/ine distributions by Jakubczyk and Respondek~\cite{JR91}, Respondek~\cite{Respondek98}, and Zhitomirskii and Respondek \cite{ZR98}.  However, intermediate ranks are important as well~-- for instance, for understanding second-order systems arising in mechanics, e.g.
\[ \ddot{x} = \sum_{i=1}^s u^i F_i(x). \]
This system can be rewritten as the following f\/irst-order system on the cotangent bund\-le~$T^*\scX$, with local coordinates $(x,p)$:
\begin{gather*}
\dot{x}   = p, \\
\dot{p}   = \sum_{i=1}^s u^i F_i(x).
\end{gather*}
The associated rank $s$ af\/f\/ine distribution fails to be strictly af\/f\/ine precisely along the submanifold $\{p=0\}$ of stationary points.  The behavior of trajectories near such points is of vital interest in control theory, so it would be useful to understand the geometry of such structures.  Elkin has introduced the notion of a ``$t$-codistribution" in \cite{Elkin98, Elkin99}, and we anticipate that this will be a~useful tool for extending our methods to af\/f\/ine distributions of this type.
\end{enumerate}

\subsection*{Acknowledgements}
This work was partially supported by NSF
grant DMS-0908456.

\pdfbookmark[1]{References}{ref}
\LastPageEnding

\end{document}